\documentclass[11pt]{amsart}
\reversemarginpar
\newcommand{\Tea}{T_{\vep,\alpha}}

\usepackage{amsmath} 
\usepackage{amssymb}

\newcommand{\nwc}{\newcommand}

\nwc{\nwt}{\newtheorem}

\nwt{prop}{Proposition}
\nwt{proposition}{Proposition}
\nwt{lem}{Lemma}
\nwt{thm}{Theorem}
\nwt{emp}{}
\nwt{cor}{Corollary}
\nwt{rem}{Remark}
\nwt{defin}{Definition} %def is already defined

\nwc{\bal}{\begin{align}}
\nwc{\be}{\begin{equation}}
\nwc{\ben}{\begin{equation*}}
\nwc{\bea}{\begin{eqnarray}}
\nwc{\beq}{\begin{eqnarray}}
\nwc{\bean}{\begin{eqnarray*}}
\nwc{\beqn}{\begin{eqnarray*}}
\nwc{\beqast}{\begin{eqnarray*}}

\nwc{\eal}{\end{align}}
\nwc{\ee}{\end{equation}}
\nwc{\een}{\end{equation*}}
\nwc{\eea}{\end{eqnarray}}
\nwc{\eeq}{\end{eqnarray}}
\nwc{\eean}{\end{eqnarray*}}
\nwc{\eeqn}{\end{eqnarray*}}
\nwc{\eeqast}{\end{eqnarray*}}

\nwc{\nn}{\nonumber}
\nwc{\bm}{\boldmath}
\nwc{\mf}{\mathbf}
\nwc{\ml}{\mathcal}

\nwc{\vep}{\varepsilon}
\nwc{\ep}{\varepsilon}
\nwc{\veps}{\varepsilon}
\nwc{\vrho}{\varrho}
\nwc{\orho}{\bar\varrho}
\nwc{\ou}{\bar u}
\nwc{\vpsi}{\varpsi}
\nwc{\lamb}{\lambda_\varepsilon}

\nwc{\eps}{\epsilon}
\nwc{\lam}{\lambda}
\nwc{\del}{\delta}

\newcommand{\commentout}[1]{}

\nwc{\la}{\langle}    %like '<' - bra
\nwc{\ra}{\rangle}    %like '>' - cat

\nwc{\IA}{\mathbb{A}} %algebraic
\nwc{\IB}{\mathbb{B}} 
\nwc{\IC}{\mathbb{C}} %complex
\nwc{\ID}{\mathbb{D}} %Dedekind
\nwc{\IE}{\mathbb{E}} %Euklides
\nwc{\IF}{\mathbb{F}} %finite field
\nwc{\IG}{\mathbb{G}} %Gauss
\nwc{\IN}{\mathbb{N}} %natural
\nwc{\IQ}{\mathbb{Q}} %rational
\nwc{\IR}{\mathbb{R}} %real 
\nwc{\IT}{\mathbb{T}} %torus
\nwc{\IZ}{\mathbb{Z}} %integers

\nwc{\m}{\mbox}
\nwc{\re}{\hbox{Re}}
\nwc{\e}{{\bf e}}
\nwc{\cG}{{\cal G}}
\nwc{\ubm}{\unboldmath}
\nwc{\mt}{\bar{t}}
\nwc{\bmxi}{\m{\bm $\xi$\ubm}}
\nwc{\bmeta}{\m{\bm $\eta$\ubm}}
\nwc{\bma}{\m{\bm $a$\ubm}}
\nwc{\bmb}{\m{\bm $b$\ubm}}
\nwc{\bmc}{\m{\bm $c$\ubm}}
\nwc{\bmd}{\m{\bm $d$\ubm}}
\nwc{\bme}{\m{\bm $e$\ubm}}
\nwc{\bmf}{\m{\bm $f$\ubm}}
\nwc{\bmg}{\m{\bm $g$\ubm}}
\nwc{\bmh}{\m{\bm $h$\ubm}}
\nwc{\bmi}{\m{\bm $i$\ubm}}
\nwc{\bmj}{\m{\bm $j$\ubm}}
\nwc{\bmk}{\m{\bm $k$\ubm}}
\nwc{\bml}{\m{\bm $l$\ubm}}
\nwc{\bmn}{\m{\bm $n$\ubm}}
\nwc{\bmo}{\m{\bm $o$\ubm}}
\nwc{\bmp}{\m{\bm $p$\ubm}}
\nwc{\bmq}{\m{\bm $q$\ubm}}
\nwc{\bmr}{\m{\bm $r$\ubm}}
\nwc{\bms}{\m{\bm $s$\ubm}}
\nwc{\bmt}{\m{\bm $t$\ubm}}
\nwc{\bmu}{\m{\bm $u$\ubm}}
\nwc{\bmv}{\m{\bm $v$\ubm}}
\nwc{\bmw}{\m{\bm $w$\ubm}}
\nwc{\bmx}{\m{\bm $x$\ubm}}
\nwc{\bmxt}{\m{\bm $x$\ubm}^\varepsilon (t)}
\nwc{\bmy}{\m{\bm $y$\ubm}}
\nwc{\bmz}{\m{\bm $z$\ubm}}

\nwc{\bmX}{\m{\bm $X$\ubm}}
\nwc{\bmB}{\m{\bm $B$\ubm}}
\nwc{\bmS}{\m{\bm $S$\ubm}}
\nwc{\bmR}{\m{\bm $R$\ubm}}
\nwc{\bmU}{\m{\bm $U$\ubm}}
\nwc{\bmE}{\m{\bm $E$\ubm}}
\nwc{\bmF}{\m{\bm $F$\ubm}}
\nwc{\bmH}{\m{\bm $H$\ubm}}
\nwc{\bmI}{\m{\bm $I$\ubm}}
\nwc{\bmP}{\m{\bm $P$\ubm}}
\nwc{\bmM}{\m{\bm $M$\ubm}}
\nwc{\bmJ}{\m{\bm $J$\ubm}}
\nwc{\bmK}{\m{\bm $K$\ubm}}
\nwc{\bmA}{\m{\bm $A$\ubm}}
\nwc{\bmD}{\m{\bm $D$\ubm}}
\nwc{\bmG}{\m{\bm $G$\ubm}}

\nwc{\bmtheta}{\m{\bm $\theta$\ubm}}
\nwc{\bmnu}{\m{\bm $\nu$\ubm}}
\nwc{\bmom}{\m{\bm $\omega$\ubm}}
\nwc{\om}{\omega}
\nwc{\Om}{\m{\bm $\Omega$\ubm}}
\nwc{\bmsigma}{\m{\bm $\sigma$\ubm}}
\nwc{\bmnabla}{\m{\bm $\nabla$\ubm}}
\nwc{\bmLambda}{\m{\bm $\Lambda$\ubm}}
\nwc{\bmlambda}{\m{\bm $\lambda$\ubm}}
\nwc{\bmtau}{\m{\bm $\tau$\ubm}}
\nwc{\bmPhi}{\m{\bm $\Phi$\ubm}}
\nwc{\bmphi}{\m{\bm $\phi$\ubm}}
\nwc{\bmPsi}{\m{\bm $\Psi$\ubm}}
\nwc{\bmpsi}{\m{\bm $\psi$\ubm}}
\nwc{\bmGamma}{\m{\bm $\Gamma$\ubm}}
\nwc{\bGamma}{\m{\bm $\Gamma$\ubm}}
\nwc{\bmgamma}{\m{\bm $\gamma$\ubm}}
\nwc{\bmQ}{\m{\bm $Q$\ubm}}

\nwc{\bE}{{\bf E}}
\nwc{\bW}{{\bf W}}
\nwc{\bF}{{\bf F}}
\nwc{\bD}{{\bf D}}
\nwc{\bJ}{{\bf J}}
\nwc{\bK}{{\bf K}}
\nwc{\bI}{{\bf I}}
\nwc{\bG}{{\bf G}}
\nwc{\bA}{{\bf A}}
\nwc{\bZ}{{\bf Z}}
\nwc{\bB}{{\bf B}}
\nwc{\bH}{{\bf H}}
\nwc{\bR}{{\bf R}}
\nwc{\ia}{{\it a}}
\nwc{\bC}{{\bf C}}
\nwc{\bx}{{\bf x}}
\nwc{\bq}{{\bf q}}
\nwc{\bfe}{{\bf e}}
\nwc{\by}{{\bf y}}
\nwc{\vy}{{\bf y}}
\nwc{\bk}{{\bf k}}
\nwc{\ba}{{\bf a}}
\nwc{\bb}{{\bf b}}
\nwc{\bz}{{\bf z}}
\nwc{\bp}{{\bf p}}
\nwc{\bS}{{\bf S}}
\nwc{\bi}{{\bf i}}
\nwc{\bw}{{\bf w}}

\newcommand{\bu}{{\bf u}}
\newcommand{\bc}{{\bf c}}
\newcommand{\bv}{{\bf v}}

\nwc{\cA}{{\cal A}}
\nwc{\vcalF}{{\cal F}}
\nwc{\cS}{{\cal S}}
\nwc{\cV}{{\cal V}}
\nwc{\cL}{{\cal L}}
\nwc{\cB}{{\cal B}}
\nwc{\cC}{{\cal C}}
\nwc{\cD}{{\cal D}}
\nwc{\cO}{{\cal O}}
\nwc{\cao}{{\cal A}^{-1}}
\nwc{\cE}{{\cal E}}
\nwc{\cf}{{\cal F}}
\nwc{\cg}{{\cal G}}
\nwc{\cH}{{\cal H}}
\nwc{\cQ}{{\cal Q}}
\nwc{\cI}{{\cal I}}
\nwc{\cJ}{{\cal J}}
\nwc{\cK}{{\cal K}}
\nwc{\cl}{{\cal L}}
\nwc{\cle}{{\cal L}^\varepsilon}
\nwc{\clu}{{\cal L}{\cal U}}
\nwc{\cm}{{\cal M}}
\nwc{\cn}{{\cal N}}
\nwc{\co}{{\cal O}}
\nwc{\cp}{{\cal P}}
\nwc{\cpt}{{\cal P}^\varepsilon_t}
\nwc{\cq}{{\cal Q}}
\nwc{\calr}{{\cal R}}
\nwc{\cs}{{\cal S}}
\nwc{\ct}{{\cal T}}
\nwc{\cu}{{\cal U}}
\nwc{\cv}{{\cal V}}
\nwc{\cw}{{\cal W}}
\nwc{\cx}{{\cal X}}
\nwc{\cy}{{\cal Y}}
\nwc{\cz}{{\cal Z}}

\nwc{\noi}{\noindent}
\nwc{\non}{\nonumber}
\nwc{\half}{\frac{1}{2}}
\nwc{\third}{\frac{1}{3}}

\nwc{\uP}{{\em \bf Proof: }}
\nwc{\uT}{\underline{Theorem:}}

\begin{document}

\title{Noise induced dissipation in 
Lebesgue-measure preserving maps on $d-$dimensional torus}

\author{Albert Fannjiang \and
Lech Wo{\l}owski}
\thanks{Department of Mathematics, 
University of California at Davis,
Davis, CA 95616,
Internet: fannjian@math.ucdavis.edu, wolowski@math.ucdavis.edu.
The research of AF is supported in part by the grant from U.S. National
Science Foundation, DMS-9971322 and UC Davis Chancellor's Fellowship}

\begin{abstract}
We consider dissipative systems resulting
from the Gaussian and $alpha$-stable noise perturbations of
measure-preserving maps on the $d$ dimensional torus. 
We study the dissipation time scale and its physical implications
as the noise level $\vep$ vanishes.

We show that nonergodic maps give rise to an $O(1/\vep)$ dissipation time
whereas ergodic toral automorphisms, including cat maps and their
$d$-dimensional generalizations, have an $O(\ln{(1/\vep)})$ dissipation 
time with a constant related to 
the minimal, {\em dimensionally averaged entropy} among the automorphism's
irreducible blocks.
Our approach reduces the calculation of the dissipation time
to a nonlinear, arithmetic optimization problem which is solved
asymptotically
by means of some fundamental theorems in theories
of convexity, Diophantine 
approximation and arithmetic progression.
We show that the same asymptotic can be reproduced by
degenerate noises as well as
mere coarse-graining.
We also discuss the implication of the dissipation time in
kinematic dynamo.

\end{abstract}
\maketitle
\section{Introduction}
\label{sI}

Irreversibility and approach to equilibrium are fundamental problems
in statistical mechanics and dynamical systems and its complete solution
is still elusive (see, e.g., \cite{Kr}). 
There are possibly many routes to irreversibility.

One view is that macroscopic systems are exceedingly difficult to isolate
from their environments for a time comparable to their dynamical time scales.
The noise as a result of interaction with environment may further trigger
irreversibility, such as approach to equilibrium,  in the systems.
The initial uncertainty involved
in preparing a physical system and the random perturbation due to measurements
as well as Gibbs' coarse-graining procedure can all be viewed
as certain noises. 
The point is that noises, intrinsic as a result of internal  
stochasticity
or extrinsic as a result of random influence from surrounding
environment, can induce effects
that would be weak or absent without noise. 

In this paper we investigate one such effect, called
dissipation, for discrete time, conservative dynamical systems 
under the influence of noise.
In particular we study the time scale, called the dissipation time,
on which the
dissipation as measured in $L^p-$norm, $1<p<\infty$,
has an {\em order one} effect even as the magnitude of noise vanishes.
Clearly the dissipation time depends on the ergodic
properties of the noiseless dynamics
as well as the noise level.

The noisy dynamical system considered in this paper can 
be viewed as a discrete generalization of the dynamics of a passive scalar 
in a 
periodic, incompressible velocity field $\mathbf{v}$
\bea
\label{SDE}
 &&d\mathbf{x}^{\vep}(t) =\mathbf{v}(\mathbf{x}^{\vep}(t))dt + \sqrt{\vep}d\mathbf{w}(t) \\
 &&\nabla \cdot \mathbf{v}(\mathbf{x})=0  \notag,
\eea
where the standard Brownian motion
$\mathbf{w}$  and the molecular diffusivity $\vep$ 
represent the stochastic perturbations as a result of
random molecular collisions (see, e.g., \cite{Ki}, \cite{F}).
The discrete-time dynamical system will be defined on 
the $d$-dimensional torus $\IT^d=\IR^d/\IZ^d$. The velocity field $\mathbf{v}$
will be replaced by arbitrary Lebesgue-measure 
preserving map $F$ defined on $\IT^{d}$ (periodicity condition).
 
In order to study the dynamics generated by $F$ it is useful to consider
its Koopman operator $U_{F}$
defined by a composition $U_{F}f:=f \circ F$, with $f$ belonging
to some Banach space of functions on $\IT^d$.
We will be mainly concerned with the standard Banach spaces $L^p(\IT^d), 1\leq p\leq \infty,$
and their subspaces $L^p_0(\IT^d)$ of functions  with zero mean $\la f\ra=0$, where $\la f\ra$ denotes
the average of $f$ w.r.t. the Lebesgue measure. In case of $L^{1}(\IT^{d})$ one can consider 
$U_{F}$ as the Frobenius-Perron operator associated with $F^{-1}$.

In the time-discrete version we
consider general $\alpha-$stable noise operator 
$G_{\vep,\alpha}:L_{0}^{2}(\IT^{d}) \mapsto L_{0}^{2}(\IT^{d})$, 
with $\alpha\in (0,1]$, defined by means 
of the Fourier transform of corresponding $\alpha$-stable noise kernel $g_{\vep,\alpha}$
\ben
G_{\vep,\alpha}f(\bx) = \int_{\IT^{d}}g_{\vep,\alpha}(\bx-\by)f(\by)d\by 
= \sum_{\bk \in \IZ^{d}}e^{-\epsilon |\bk|^{2\alpha}}\hat{f}(\bk)\bfe_{\bk}(\bx), 
\een
where 
\ben
\label{kernel}
g_{\vep,\alpha}(\bx):=\sum_{\bk \in \IZ^{d}}e^{-\epsilon |\bk|^{2\alpha}}\bfe_{\bk}(\bx),
\een
with $\bfe_{\bk}(\bx):=e^{2 \pi i \bk \cdot \bx}$, $\bk \in \IZ^{d}$.
Here, just like in (\ref{SDE}), $\vep>0$ represents the level of the noise. 
Putting $\alpha=1$ one recovers standard heat kernel.

The operator $T_{\vep,\alpha}$ on $L_{0}^2(\IT^{d})$ generating the noise-perturbed dynamical 
system considered in this paper is thus given by
\bea
\label{T}
T_{\vep,\alpha}f:=G_{\vep,\alpha}U_{F}f=g_{\vep,\alpha}*(f \circ F).
\eea   
Simple computations yield
\bea
\label{contr}
\|T_{\vep,\alpha}^{n}\|=\|G_{\vep,\alpha}U_{F}...G_{\vep,\alpha}U_{F}\| \leq 
\|G_{\vep,\alpha}\|^{n}=e^{-\vep n}.
\eea
Here and throughout the paper $\|\cdot\|$ denotes the standard $L_0^{2}$-norm
or the corresponding operator 
norm (any other norm will be equipped with suitable index).

We define the dissipation time as the time
on which
the contraction (\ref{contr}) becomes of order one:
\bea
\label{tdiss1}
n_{diss}:=\min\{n \in \IZ_{+}:\|T_{\vep,\alpha}^{n}\| < 1/e\}.
\eea
Hence the dissipation time is a function of $\ep,\alpha$ as well as the
underlying dynamics. The choice of the threshold $e^{-1}$ in the definition is
a convenient one and for the purpose of the
paper can be any positive number less than one (see Proposition~\ref{prop1}).
The fact that for all $\vep>0$, $\|T_{\vep,\alpha}^{n}\|$ is monotonically decreasing
ensures that $n_{diss}$ is well defined.
By contrast, when $\vep=0$,
the fine-grained Boltzmann-Gibbs entropy as
well as $L^p$-norm of the initial state remains
constant in the course of evolution.
In other words this is a ``dissipation'' effect and
hence the term ``dissipation time''.
 On the dissipation time scale the system is, in a sense, ``half way''
 through
 its irreversible route to the equilibrium state. The dissipation time
 provides a measure of the instability of the dynamics w.r.t.
 the stochastic perturbations which result in the ``aging'' of
the system toward
the final state. 

The main purpose of this paper is to investigate  
the asymptotics of the dissipation time as $\ep$ tends to zero.
%The asymptotic behavior of $n_{diss}$ provides, among other things,
%a sharp lower bound
%for the time scale on which the noisy dynamics can be approximated by
%a Brownian motion regardless whether the unperturbed
%system is ergodic or not.
Due to the non-normality of the operator
 $T_{\vep,\alpha}$, the dissipation time can not be determined
from its spectral radius.
 Indeed, as we will see below, the operator $T_{\vep,\alpha}$ corresponding
  to any ergodic toral automorphism $F$ is quasi-nilpotent for any $\vep>0$
   (thus the time scale estimated from the spectral radius is infinite)
 whereas the dissipation time is of the order $\ln{(1/\ep)}$.

To briefly describe the main results 
we pause to note the following asymptotic notation.
Given two sequences $a_{\vep},b_{\vep}$ indexed by the parameter $\vep>0$
we write
\bean
a_{\vep} &\lesssim& b_{\vep} \text{, if } \limsup_{\vep
\rightarrow
0}\frac{a_{\vep}}{b_{\vep}} < \infty \\
a_{\vep} &\approx& b_{\vep} \text{, if } \lim_{\vep \rightarrow 0
 }\frac{a_{\vep}}{b_{\vep}}=1
\eean
 Moreover we write $a_{\vep} \sim b_{\vep}$, if both 
$a_{\vep} \lesssim b_{\vep}$ and $b_{\vep} \lesssim a_{\vep}$ hold simultaneously.

The first result is that
the dissipation time $n_{diss}\sim 1/\vep$ for
nonergodic or, more general, non-weakly-mixing maps (cf. Theorem \ref{nonergodic} and its
corollaries, Section \ref{sDTA}),
which is also the longest possible time scale for dissipation in view of
(\ref{contr}) (see also Lemma \ref{tdlub}, Section  \ref{sDCal}).
In other word, such systems are most stable w.r.t. stochastic perturbations.

The main aim, however, is to investigate the cases
in which the dissipation is much faster due to rapid mixing in $F$.
We
show (Theorem 2, Section \ref{sDTA})
that, for a toral automorphism $F$, $n_{diss} \sim \log(1/\vep)$
if and only if the map $F$ is ergodic (which in this case is also an
Anosov diffeomorphism).
In particular our results hold for all classical cat maps (hyperbolic
automorphisms of $2$-torus) and their $d$-dimensional generalizations. 
In addition, we provide a general lower bound for the constant of the
logarithmic asymptotics (Theorem 2). We further show
that the lower bound is achieved for 
{\em diagonalizable} automorphisms, namely
\be
\label{open}
n_{diss}\approx (2\alpha \hat{h}(F))^{-1} \ln(1/\vep)
\ee
where $\hat{h}(F)$ denotes the minimal, dimensionally
averaged entropy among
$F$'s irreducible blocks (Theorem 3, Section \ref{sDTA}).
Dimensionally averaged 
entropy for each irreducible sub-block of the
toral automorphism is 
the Kolmogorov-Sinai (KS) entropy per dimension of an irreducible factor
of the whole map.
%It relates to the
%factorization of the Hilbert space of the square integrable functions 
%into tensor product of
%invariant subspaces w.r.t. the underlying noiseless dynamics.  
%Essentially all mixing Anosov diffeomorphisms should have 
%the $\log{1/\ep}$ dissipation time but it
%is not immediately clear
%what the constant
%should be.

Our method involves  solving asymptotically a
quadratic arithmetic optimization (i.e. quadratic integer programming)
problem by obtaining sharp upper and
lower bounds using number theoretical
tools including multidimensional Diophantine approximation theorems (Schmidt's
subspace theorem), Minkowski's theorem on linear forms and Van der Waerden's theorem 
on arithmetic progressions. This is done in Section~\ref{sA}. 

In Section~\ref{coarse} we show that the same result (\ref{open})
holds when the noise is replaced by  coarse-graining
the initial and terminal states.
This is reminiscent of the well-know results
of statistical stability in the literature, namely,
the Bernoulli systems are stable under the
sufficiently small intrinsic random perturbation 
in the rough sense that
the perturbed system is close to the direct product of
the unperturbed one and some auxiliary viewer system
(see \cite{Ki}, \cite{OW}).
In other words, for those systems, the process that results
from small intrinsic random perturbation can be
reproduced exactly by looking at the unperturbed system
through a viewer that distorts randomly but slightly.
In spite of the above 
the asymptotic (\ref{open}) indicates
the perturbed system is irreversibly far from the unperturbed one
even on the relatively short dissipation time scale.
From this perspective, such  a system is
statistically unstable.

In Section \ref{degnoise} we consider a class
of highly degenerate 
noises and show that the same conclusions about the dissipation
hold if the degenerate noises satisfy an additional generic condition.   

In Section \ref{dynamo} we consider the relation
between the dissipation time and some characteristic
time scales relevant to kinematic dynamo. 
We show that fast dissipation generally inhibits dynamo action.
When there is no fast dynamo action
the noisy push-forward map 
dissipates the magnetic field energy on the dissipation time  scale.
However, the magnetic field energy can still grow to
relatively large magnitude as inverse power-law of
the small noise with the exponent proportional to
the ratio of the logarithmic spectral radius of the toral automorphism to
the minimal, dimensionally averaged entropy among
the automorphism's irreducible blocks (cf. (\ref{anti})).

The notion of dissipation time has a natural bearing
on the problems of quantum chaos with noise.
The family of symplectic toral automorphisms constitute 
important examples of quantizable
chaotic dynamics on compact manifolds for which various
quantization procedures have been intensively studied (see, for example,
\cite{RO} and \cite{Mez}). 
We will address the issue of decoherence time for
quantized symplectic toral automorphisms with noise in a  forthcoming paper. 

The organization of the rest of the paper is as follows.
In Section \ref{sD} we develop the general 
theory of dissipation time and its relation to the Boltzmann-Gibbs entropy.
We also formulate the dissipation time calculation for
total automorphisms as an arithmetic minimization problem and state
the main results.
In Appendix \ref{aAT} we generalize the dissipation time asymptotic
result
to the affine transformations. 
The proofs of some elementary facts are presented in Appendix \ref{aP} 
for the sake of completeness.

\section{Dissipation time}
\label{sD}

In its general form the dissipation time $n_{diss}(p)$ 
can be defined in terms of the norm $\|\cdot\|_{p,0}$
 on the space $L^p_0(\IT^d)$ w.r.t. a threshold $\eta\in(0,1)$
\bea
\label{tdiss2}
n_{diss}(p,\eta):=\min\{n \in \IZ_{+}:\|T_{\vep,\alpha}^{n}\|_{p,0}< \eta\}, \quad 1\leq p\leq\infty.
\eea
First we show that the value of the threshold $\eta$ in (\ref{tdiss2}) does not affect 
the order of divergence of $n_{diss}(p,\eta)$, as $\vep$ tends to zero.
\begin{prop}
\label{prop1}
For any $0<\tilde{\eta},\eta<1$,  $n_{diss}(p,\tilde{\eta}) \sim n_{diss}(p,\eta)$.
\end{prop}
\textbf{Proof.}
Assume $0<\tilde{\eta}<\eta<1$. 
Obviously $n_{diss}(p,\tilde{\eta}) \geq n_{diss}(p,\eta)$. On the other hand let $k$ be a 
positive integer such that $\eta^{k}<\tilde{\eta}$. Then
\be
 \|T_{\vep,\alpha}^{n_{diss}(p,\eta)}\|_{p,0}<\eta \Rightarrow \|T_{\vep,\alpha}^{kn_{diss}(p)}\|_{p,0}
<\eta^{k}<\tilde{\eta}.
\ee
Hence $kn_{diss}(p,\eta) \geq n_{diss}(p,\tilde{\eta})$, which implies 
$n_{diss}(p,\eta) \sim n_{diss}(p,\tilde{\eta})$.$\qquad \blacksquare$

Following the argument of \cite{Ros} one can use the Riesz convexity theorem to establish
also the asymptotic equivalence of the $n_{diss}(p)$, for all $1<p<\infty$.
\begin{prop}
\label{pdiss}

i) For any $1<q,p<\infty$,  $n_{diss}(q) \sim n_{diss}(p)$.

ii) For any $1<p<\infty$, $n_{diss}(p) \lesssim n_{diss}(1)$ and $n_{diss}(p) \lesssim n_{diss}(\infty)$.
\end{prop}

The details of the proof can be found in Appendix \ref{aP}.

Our particular choice of the exponent $p=2$ and threshold $\eta=e^{-1}$ in (\ref{tdiss1}) 
is computationally convenient and will be used throughout the paper. 
We will use the convention that $n_{diss}(p)=n_{diss}(p,e^{-1})$.

We say that operator $T_{\vep,\alpha}$ or associated with it measure preserving map $F$ has 
a {\em simple (slow)} dissipation time when $n_{diss} \sim 1/\vep$
and that it has a {\em logarithmic (fast)} dissipation time when $n_{diss} \sim \ln(1/\vep)$.

In the particular case of fast dissipation, with a logarithmic dissipation time, 
in order to estimate precisely the rate of dissipation, one needs to determine the value of the 
{\em dissipation rate constant} $R_{diss}$, defined as 
\bea
\label{Rdiss}
  R_{diss}:=\lim_{\vep \rightarrow 0} \frac{n_{diss}}{\ln(1/\vep)}.
\eea    
Similarly in case of simple dissipation time the dissipation rate 
constant can be defined as
\bea
\label{Rdissnon}
  R_{diss}:=\lim_{\vep \rightarrow 0} \vep n_{diss}.
\eea  

\subsection{Dissipation time and Boltzmann-Gibbs entropy}
\label{sDB-G}
In this section we briefly discuss the connection between
dissipation time and Boltzmann-Gibbs entropy.

First we note that on the scale of $n_{diss}$ the Boltzmann-Gibbs entropy
approaches the maximal equilibrium value (i.e. $0$) as can
be seen from the following simple estimate \cite{LM}.
Let us first restrict considerations to bounded initial states, i.e., $f\geq 0, f\in L^\infty$
and $\|f\|_{1}=1$. Let 
\[
\eta(u)=\begin{cases}
-u\ln{u},& u>0\\
0,& u=0
\end{cases}
\]
and let $D_n=\{\bx\in \IT^d: 1\leq T^n_{\vep,\alpha}f\}$. On one hand, we have
\bea
&&\left|\int_{D_n}\eta(T^n_{\vep,\alpha}f(\bx))d\bx\right| \label{fi}\nonumber \\
&\leq&\int_{D_n}\left|\int_{1}^{T^n_{\vep,\alpha}f(\bx)}\frac{d\eta(u)}{du}du\right|
d\bx \notag \\
&\leq&\sup_{1\leq u\leq \|\Tea^n f\|_\infty}(1+\ln{u})\int_{D_n}
|\Tea^n f(\bx)-1|d\bx \notag \\
&\leq& (1+\ln{\|\Tea^n f\|_\infty}) \|\Tea^nf-1\|_1 \notag \\
&\leq&(1+\ln{\| f\|_\infty})\|\Tea^nf-1\|_1 \label{li}.
\eea
On the other hand, we have
\[
0 \geq \int_{\IT^d}\eta(\Tea^n f(\bx))d\bx 
\geq \int_{D_n}\eta(\Tea^n f(\bx))d\bx.
\]
In view of the inclusion relation: $L^\infty(\IT^d)\subset L^2(\IT^d)
\subset L^1(\IT^d)$, we then obtain that for $n\gg n_{diss}$ 
\[
\sup_{f\geq 0,\|f\|_\infty\leq c}\left|\int_{\IT^d}\eta(\Tea^n f(\bx))d\bx
\right|
\stackrel{\vep \downarrow 0}{\longrightarrow} 0,
\quad \forall c>0.
\]

For unbounded initial states, we note that, by Young's inequality,
\[
 \|\Tea^n f\|_\infty \leq  \|\Tea f\|_\infty \leq 
 \|g_{\vep,\alpha}\|_\infty \|f\|_{1}=\|g_{\vep,\alpha}\|_\infty
\] 
from which we have, instead of (\ref{li}), the following estimate 
\[
\left|\int_{D_n}\eta(T^n_{\vep,\alpha}f(\bx))d\bx\right| \leq
(1+\ln{\|g_{\vep,\alpha}\|_\infty})\|\Tea^nf-1\|_1.
\] 
where
\[
 \ln\|g_{\vep,\alpha}\|_\infty
 \sim \ln(1/\vep).
\]
Therefore for sufficiently fast diverging $n\gg n_{diss}(1)$ such 
that
\be
\label{entropy2}
 \ln(1/\vep)\|\Tea^n(f-1)\|_{1,0}
 \stackrel{\vep \downarrow 0}{\longrightarrow} 0
\ee 
one obtains
\[
\sup_{f\geq 0,\|f\|_1=1}\left|\int_{\IT^d}\eta(\Tea^n f(\bx))d\bx
\right|
\stackrel{\vep \downarrow 0}{\longrightarrow} 0.
\]
The condition (\ref{entropy2}) typically results in
a slightly longer time scale than $n_{diss}(1)$.

On the other hand, we  can bound the $L_1$ distance between
the probability density function $f$ and the Lebesgue measure
by their relative entropy via Csisz\'{a}r's inequality \cite{Cs}
\[
\int_{\IT^d} 
|f(\bx)-g(\bx)|d\bx \leq \sqrt{2\int_{\IT^d} f(\bx)\ln{(f(\bx)/g(\bx))}d\bx}
\]
with $g(\bx)=1$. We see immediately that the decay rate
of 
\[
\sup_{f\geq 0,\|f\|_1=1}\left|\int_{\IT^d}\eta(\Tea^n f(\bx))d\bx
\right|
\]
provides an estimate for $n_{diss}(1)$ and, consequently, for
$n_{diss}(p), p\in (1,\infty)$.

\subsection{Calculating the dissipation time}
\label{sDCal}
For greater generality and transparency of arguments we consider, in this section, a slightly more general family of operators $T_{\vep,\alpha}$ defined, as previously, by the first equality in (\ref{T}), but with  
arbitrary unitary or isometric (not necessary Koopman) operator $U$ (and hence
in these cases we drop the subscript $F$). 

\begin{lem}
\label{tdlub}
For any isometric operator $U$, the dissipation time of $T_{\vep,\alpha}$ satisfies following constraints 
\bea
\label{tdisslub}
\|R(1;T_{\vep,\alpha})\| \lesssim n_{diss} \lesssim 1/\vep,
\eea 
where $R(1;T_{\vep,\alpha})$ denotes the resolvent of $T_{\vep,\alpha}$ at 1.
\end{lem}
\textbf{Proof.}
In view of (\ref{tdiss1}) and (\ref{contr}), for $n=n_{diss}$ one has
\bean
e^{-1} \leq \|T_{\vep,\alpha}^{(n_{diss}-1)}\| \leq e^{-\vep (n_{diss}-1)},
\eean
which clearly implies the second estimate of (\ref{tdisslub}).
In order to prove the other inequality we proceed as follows.
\bean
 &&\|R(1;T_{\vep,\alpha})\|=
 \left\|\sum_{n=0}^{\infty}T_{\vep,\alpha}^{n}\right\|=
 \left\|\sum_{n=0}^{n_{0}-1}T_{\vep,\alpha}^{n} + T_{\vep,\alpha}^{n_{0}}\sum_{n=0}^{\infty}T_{\vep,\alpha}^{n}\right\|\\
 &\leq&\sum_{n=0}^{n_{0}-1}\|T_{\vep,\alpha}^{n}\| + \|T_{\vep,\alpha}^{n_{0}}\|\left\|\sum_{n=0}^{\infty}T_{\vep,\alpha}^{n}\right\| \leq
 n_{0} + \|T_{\vep,\alpha}^{n_{0}}\|\|R(1;T_{\vep,\alpha})\|.
\eean
Hence taking in the above inequality $n_{0}=n_{diss}$ one gets
\bean
 \|R(1;T_{\vep,\alpha})\|(1-e^{-1}) \leq
 \|R(1;T_{\vep,\alpha})\|(1-\|T_{\vep,\alpha}^{n_{diss}}\|) 
 \leq n_{diss},
\eean
which gives the first estimate of (\ref{tdisslub}). $\qquad \blacksquare$

The above lemma provides an absolute upper bound for dissipation time. Taking $F=I$ one easily  
finds that this bound is best possible in general. The lower bound is useful in the case when one can
estimate from below the norm of the resolvent (see proof of Theorem \ref{nonergodic}).  

\begin{thm}
\label{nonergodic}
If $U$ acting on $L^{2}_{0}(\IT^{d})$ possesses nonempty pure point spectrum 
and at least one of its eigenfunctions belongs to $H^{2\alpha}(\IT^{d})$, then $T_{\vep,\alpha}$ 
has simple dissipation time.
\end{thm}
\textbf{Proof}.
In view of Lemma~\ref{tdlub} it is enough to find a lower bound for the norm of the 
resolvent $R(1;T_{\vep,\alpha})$. Let $h \in H^{2\alpha}$ be one of the eigenfunctions of $U$. 
Since $U$ is isometric we have 
\bean
Uh=e^{i\phi}h.
\eean
We first assume that $\phi=0$. 
Since $1 \not \in \sigma (T_{\vep,\alpha})$, $I-T_{\vep,\alpha}$ is a homeomorphism and hence
\bean
 \|R(1;T_{\vep,\alpha})\|=
 \sup_{f\in L^{2}_{0}}\frac{\|(I-T_{\vep,\alpha})^{-1}f\|}{\|f\|}=
 \sup_{f\in L^{2}_{0}}\frac{\|f\|}{\|(I-T_{\vep,\alpha})f\|} \geq
 \frac{\|h\|}{\|(I-T_{\vep,\alpha})h\|}.
\eean
Now expressing $h$ in the Fourier series we get
\bean
 \|(I-T_{\vep,\alpha})h\|^{2}=
 \sum_{0 \not =\bk\in \IZ^{d}}|\hat{h}(\bk)|^{2}\left|1-e^{-\vep|\bk|^{2\alpha}}\right|^{2} \leq
 \sum_{0 \not =\bk\in \IZ^{d}}\left(\vep|\hat{h}(\bk)||\bk|^{2\alpha}\right)^{2}=
 \vep^{2} \|h\|_{H^{2\alpha}}^{2}.
\eean
Hence
\bean
 \|R(1;T_{\vep,\alpha})\| \geq \frac{\|h\|}{\vep\|h\|_{H^{2\alpha}}}=: \frac{C}{\vep},
\eean
Thus in view of (\ref{tdisslub}) and above calculations 
\bean
1/\vep \lesssim \|R(1;T_{\vep,\alpha})\| \lesssim n_{diss} \lesssim 1/\vep,
\eean
which ends the proof in the case $\phi=0$.

If $\phi \not=0$, we put
 $\hat{U}= e^{-i\phi}U$,
which implies 
$\hat{U}h=h.$

The proof is completed by applying the above reasoning to operator $\hat{T}_{\vep,\alpha}=G_{\vep,\alpha}\hat{U}$ 
and observing that the dissipation times for $T_{\vep,\alpha}$ and $\hat{T}_{\vep,\alpha}$ are identical. 
$\qquad \blacksquare$

When $U$ is a Koopman operator associated with a map $F$, then the property that 
$U$ considered on $L^{2}_{0}(\IT^{d})$ possesses nonempty pure point spectrum is equivalent to the 
fact that $F$ is not weakly mixing (see \cite{CoFoSi}). Thus we have
\begin{cor}
If $F$ is not weakly mixing and its Koopman operator possesses $H^{2\alpha}$ eigenfunction 
 in $L^2_0(\IT^{d})$, then $T_{\vep,\alpha}$ has simple dissipation time.
\end{cor}
Another immediate consequence is
\begin{cor}
\label{ne}
If $F$ is not ergodic and its nontrivial invariant measure possesses $H^{2\alpha}$ density function, then 
$T_{\vep,\alpha}$ has simple dissipation time.
\end{cor}

A typical example of ergodic but not weakly mixing transformations for which the above corollary applies
is the family of 'irrational' shifts on $\IT^{d}$ i.e. maps $F\bx=\bx+\mf{c}$ on $\IT^{d}$, where $\mf{c}=(c_{1},..,c_{d})$ is a constant vector such that the numbers $1,c_{1},..,c_{d}$ are linearly independent over rationals. 
More general and less trivial examples of ergodic maps giving rise to a simple dissipation time will be 
discussed in Appendix A (cf. Remark \ref{remerg}). 

In general the problem of computing the dissipation time is rather complicated. In some cases 
it can be reformulated as an asymptotic optimization problem. 
To see it, one can represent the action of a given unitary operator $U$ in the Fourier basis
\bean
\label{uk}
 U\bfe_{\mf{\bk}}=\sum_{0 \not = \bk' \in \IZ^{d}}u_{\bk,\bk'}\bfe_{\bk'},
\eean
where for each $\bk$
\bea
\label{ul2}
 \sum_{0 \not = \bk' \in \IZ^{d}}|u_{\bk,\bk'}|^{2}=1.
\eea
Next we introduce the notation
\bean
\ml{U}_{n}(\bk_0,\bk_{n})&=&\sum_{0 \not = \bk_{1},...,\bk_{n-1} \in \IZ^{d}}
u_{\bk_{0},\bk_{1}}...u_{\bk_{n-1},\bk_{n}} e^{-\vep \sum_{l=1}^{n}|\bk_{l}|^{2\alpha}} \\
\ml{S}_n(\bk_{n})&=&\{\bk_{0} \in \IZ^{d}\backslash \{0\} : \ml{U}_{n}(\bk_0,\bk_{n})\not = 0\}.
\eean
Then for any $f\in L_{0}^{2}(\IT^{d})$  we have
\bea
\nonumber
\|T_{\vep,\alpha}^{n}f\|^{2} &=& 
\left\|\sum_{0 \not = \bk_{0} \in \IZ^{d}}\hat{f}(\bk_{0})T_{\vep,\alpha}^{n} 
\bfe_{\bk_{0}}\right\|^{2}= 
\left\|\sum_{0 \not = \bk_{0} \in \IZ^{d}}\hat{f}(\bk_{0})\sum_{0 \not = \bk_{n} \in \IZ^{d}} 
\ml{U}_{n}(\bk_0,\bk_{n})\bfe_{\bk_{n}}\right\|^{2} \\
\label{Tnorm1}
&=&\sum_{0 \not = \bk_{n} \in \IZ^{d}}\left| \sum_{0 \not = \bk_{0} \in \IZ^{d}}\hat{f}(\bk_{0})
 \ml{U}_{n}(\bk_0,\bk_{n})\right|^{2}=
\sum_{0 \not = \bk_{n} \in \IZ^{d}}\left| \sum_{\bk_{0} \in \ml{S}_{n}(\bk_{n})}\hat{f}(\bk_{0})
 \ml{U}_{n}(\bk_0, \bk_{n})\right|^{2}.
\eea
The following general upper bound for $\|T_{\vep,\alpha}^{n}f\|$
holds.
\begin{lem}
\label{gub}
For any $f\in L_{0}^{2}(\IT^{d})$,
\bea
\label{uppernorm}
\|T_{\vep,\alpha}^{n}f\|^{2} &\leq&
\sum_{0 \not = \bk_{n} \in \IZ^{d}}
\sum_{\bk_{0} \in \ml{S}_{n}(\bk_{n})}
|\hat{f}(\bk_{0})|^{2}
\sum_{\bk_{0} \in \ml{S}_{n}(\bk_{n})}
|\ml{U}_n(\bk_{0},\bk_{n})|^{2}.
\eea
\end{lem}
For the proof we refer the reader to Appendix \ref{aP}.

When $u_{\bk,\bk'}$ is a Kronecker's delta function 
\bea
\label{linear}
u_{\bk,\bk'}=\del_{A\bk,\bk'},
\eea
where $A:\IZ^{d}\mapsto \IZ^{d}$ is a linear surjective map,
the upper bound (\ref{uppernorm}) can be used to obtain
an identity for $\|T_{\vep,\alpha}^{n}\|$. 
First observe that 
\bean
\label{linear1}
 \ml{U}_n(\bk_{0},\bk_{n})=e^{-\vep \sum_{l=1}^{n}|A^{l}\bk_{0}|^{2\alpha}}\del_{A^{n}\bk_{0},\bk_{n}}
\eean
and hence (\ref{uppernorm}) becomes
\bean
\|T_{\vep,\alpha}^{n}f\|^{2} &\leq&
\sum_{0 \not = \bk_{0} \in \IZ^{d}}
|\hat{f}(\bk_{0})|^{2}e^{-2\vep \sum_{l=1}^{n}|A^{l}\bk_{0}|^{2\alpha}} \leq
\|f\|^{2}\max_{0 \not=\bk \in \IZ^{d}}
e^{-2\vep \sum_{l=1}^{n}|A^{l}\bk|^{2\alpha}}.
\eean
On the other hand for any nonzero $\bk\in \IZ^{d}$, one can take in (\ref{Tnorm1}) $f=\bfe_{\bk}$   
and get
\bean
\|T_{\vep,\alpha}^{n}f\|^{2}= e^{-2\vep \sum_{l=1}^{n}|A^{l}\bk|^{2\alpha}}
\eean
and therefore
\bea
\label{Tnorm3}
\|T_{\vep,\alpha}^{n}\| &=& \max_{0 \not=\bk \in \IZ^{d}} e^{-\vep \sum_{l=1}^{n}|A^{l}\bk|^{2\alpha}} = 
e^{-\vep \min_{0 \not=\bk \in \IZ^{d}} \sum_{l=1}^{n}|A^{l}\bk|^{2\alpha}}. 
\eea

Let us now determine the class of maps $F$ such that the corresponding Koopman operator 
$U_{F}$ satisfies (\ref{linear}). The relation (\ref{linear}) implies 
\bean
U_{F}\bfe_{\bk}=\bfe_{A\bk}=e^{2\pi i \la A\bk,\bx\ra}.
\eean 
On the other hand
\bean
U_{F}\bfe_{\bk}(\bx)=\bfe_{\bk}(F\bx)=e^{2\pi i \la \bk,F\bx\ra}.
\eean
Thus 
\[
\la \bk,F\bx\ra=\la A\bk,\bx\ra \,\,\hbox{mod}\,\, 1,\quad\forall \bx \in \IR^d,\,\, \bk\in \IZ^d,
\] that is, $A$ is linear and $A^\dagger$ equals the lifting of $F$
from $\IT^{d}$ onto $\IR^{d}$. 
Moreover,  the matrix $A$ has integer entries and  determinant equal to $\pm 1$, i.e., $A$ (and $F$)
is a toral automorphism.
Hence, for toral automorphisms, the calculation of the dissipation time reduces to
 the following nonlinear, asymptotic (large $n$) 
arithmetic minimization problem 
\bea 
\label{artminA}
\min_{0 \not=\bk \in \IZ^{d}} \sum_{l=1}^{n}|A^{l}\bk|^{2\alpha}.
\eea
 We will show in Section \ref{sA} that for any ergodic toral automorphism this 
minimum value grows geometrically in $n$ with the base related to the dimensionally-averaged
KS-entropy of the total automorphism.

\subsection{Dissipation time of toral automorphisms} 
\label{sDTA}
It is well known (see \cite{AP}) that (the lifting map corresponding to) any 
toral homeomorphism $H:\IT^{d} \mapsto \IT^{d}$ can be decomposed
into three parts $H=L+P+c$, where $L$, the linear part, is an element of $SL(d,\IZ)$ - the set of all 
matrices with integer entries and determinant equal to $\pm 1$, $P$ is periodic i.e. 
$P(\bx+\bv)=P(\bx)$ for any $\bv \in \IZ^{d}$, and $c$ is a constant shift vector.

Every algebraic and measurable automorphism of the torus is continuous.
Each continuous toral automorphism is a homeomorphism with zero periodic 
and constant parts and hence can be identified with an element of $SL(d,\IZ)$. 
And vice versa, each element of $SL(d,\IZ)$ uniquely determines a measurable, algebraic 
toral automorphism.
Thus from now on the term {\em toral automorphism}
will simply be reserved for elements of $SL(d,\IZ)$.
We recall here that all Anosov diffeomorphisms on $\IT^d$ are topologically conjugate 
to the toral automorphisms (\cite{Fr}, \cite{Man}).

Below we summarize some ergodic properties of toral automorphisms 
(cf. \cite{Katok} p. 160, \cite{Katznelson} and \cite{Arnold}).
\begin{prop}
\label{ergodicTA}
Let $F$ be a toral automorphism. The following statements are equivalent

a) no root of unity is an eigenvalue of $F$.

b) $F$ is ergodic.

c) $F$ is mixing.

d) $F$ is a K-system.

e) $F$ is a Bernoulli system.
\end{prop}

In the sequel we will use the following result (cf. \cite{Yuz}).
\begin{prop} The entropy $h(F)$ of any toral endomorphism $F$ is computed by the formula
\bea
\label{entropy}
 h(F)=\sum_{|\lam_{j}| \geq 1}\ln{|\lam_{j}|},
\eea 
where $\lam_{j}$ denote the eigenvalues of $A$.
\end{prop}
>From the formula (\ref{entropy}) one immediately sees that a toral automorphism has 
zero entropy if and only if
all its eigenvalues are of modulus 1. In fact much stronger result holds.
\begin{prop}
\label{zeroentropy}
 A toral automorphism has zero entropy if and only if
 all its eigenvalues are roots of unity. 
 In particular all ergodic toral automorphisms have positive entropy.
\end{prop}

Given any toral automorphism $F$ we denote by $P$ its characteristic polynomial and by 
$\{P_{1},...,P_{s}\}$ the complete set of its distinct irreducible (over $\IQ$) factors. 
Let $d_{j}$ denote the degree of polynomial $P_{j}$ and $h_{j}$ the   
KS-entropy of any toral automorphism with the characteristic polynomial $P_{j}$. 
For each $P_{j}$ we define its dimensionally averaged KS-entropy as 
\bea
 \hat{h}_{j}=\frac{h_{j}}{d_{j}}.
\eea 
For the whole matrix $F$ we define its minimal dimensionally averaged entropy (denoted $\hat{h}(F)$) as
\bean
 \hat{h}(F)=\min_{j=1,...,s}\hat{h}_{j} 
\eean

Now we state two main theorems of the present paper.
\begin{thm}
\label{thm2}
Let $F$ be any toral automorphism, $U_{F}$ the Koopman operator
associated with $F$, $G_{\vep,\alpha}$ $\alpha$-stable noise 
operator and $T_{\vep,\alpha}=G_{\vep,\alpha}U_{F}$. Then
 
i) $T_{\vep,\alpha}$ has simple dissipation time 
if and only if $F$ is not ergodic.

ii) $T_{\vep,\alpha}$ has logarithmic dissipation time 
if and only if $F$ is ergodic. 

iii) If $T_{\vep,\alpha}$ has logarithmic dissipation time then the dissipation rate constant 
     satisfies the following constraint
\bean
\frac{1}{2 \alpha \hat{h}(F)} \leq R_{diss} \leq \frac{1}{2 \alpha \tilde{h}(F)},
\eean
where $\tilde{h}(F)$ is a positive constant satisfying $\tilde{h}(F)\leq\hat{h}(F)$. 
\end{thm}

Part i) of the above theorem follows immediately from
Theorem \ref{nonergodic}. For details of a simple proof we refer to 
appendix \ref{aP}. 

The natural question arises, whether the lower bound for the dissipation rate constant 
given in the above theorem is best possible. 
The next theorem and its corollary provides a strong argument in favor of this conjecture.

\begin{thm}
\label{thm3}
If $F$ is ergodic and diagonalizable then 
\bean
 n_{diss}\approx \frac{1}{2 \alpha \hat{h}(F)}\ln(1/\vep).
\eean
That is, the dissipation rate constant of 
$T_{\vep,\alpha}$ is given by
\bean
 R_{diss}=\frac{1}{2 \alpha \hat{h}(F)}.
\eean
\end{thm}

The proof of parts ii) and iii) of Theorem \ref{thm2} and of 
Theorem \ref{thm3} constitute the most important part of this 
work and will be presented in Section \ref{sAT} after 
necessary tools are developed. 

We end this section with the 
the results for two and three dimensional
tori. 
Ergodicity of two dimensional toral automorphisms is equivalent
to hyperbolicity. Two dimensional hyperbolic toral automorphisms are 
often referred to as the \emph{cat maps}.
  
Using Corollary \ref{23dim} and applying Theorem \ref{thm3} to
two and three dimensions one gets the following
\begin{cor}
\label{cat}
 Let $F$ be any ergodic, two or three dimensional toral automorphism. Then 
\bean
 n_{diss}\approx \frac{1}{2 \alpha \hat{h}(F)}\ln(1/\vep),
\eean
\end{cor}

\section{Asymptotic arithmetic minimization problem}
\label{sA}
In this section we find the asymptotics, as n goes to infinity, of the following 
quadratic arithmetic minimization problem
\bea 
\label{artmin0}
\min_{0 \not=\bk \in \IZ^{d}} \sum_{l=1}^{n}|A^{l}\bk|^{2\alpha},
\eea
where $A\in SL(d,\IZ)$.
When $A$ is not ergodic the asymptotics of (\ref{artmin0}) is 
of the order $O(n)$. For the rest of the paper we will consider
only the ergodic case. 
For $d=2$ the problem (\ref{artmin0}) can be solved easily
as follows. 
Consider first the case that $A$ is symmetric and $\alpha=1$.
>From $det(A)=1$ we see that eigenvalues are $\lam, \lam^{-1}$ with
$|\lam|>1$.
We have
\bean
\min_{0 \not =\bk \in \IZ^{d}} \sum_{l=1}^{2n+1}|A^{l}\bk|^{2}&=&
\min_{0 \not=\bk \in \IZ^{d}} \sum_{l=-n}^{n}|A^{l}\bk|^{2}\\
&=&\min_{0 \not=\bk=\bk_{1}+\bk_{2} \in \IZ^{d}}
\left( |\bk|^{2} + \sum_{l=1}^{n}|\lam|^{2l}|\bk_{1}|^{2}+|\lam|^{-2l}|\bk_{2}|^{2}
+\sum_{l=1}^{n}|\lam|^{-2l}|\bk_{1}|^{2}+|\lam|^{2l}|\bk_{2}|^{2}\right)\\
&=&
\min_{0 \not=\bk \in \IZ^{d}}\sum_{l=-n}^{n}|\lam|^{2l}|\bk|^{2}
=\sum_{l=-n}^{n}|\lam|^{2l}.
\eean
Hence there exist constants $C_{1}$ and $C_{2}$ such that 
\bean
C_{1}e^{h(A)n} \leq \min_{0 \not=\bk \in \IZ^{d}} \sum_{l=1}^{n}|A^{l}\bk|^{2\alpha}
\leq C_{2}e^{h(A)n}.
\eean
where $h(A)$ denotes the KS-entropy of $A$.
The estimates for the general case of non-symmetric $A$ and $\alpha \neq 1$
are similar. 

In higher dimensions, the solution to (\ref{artmin0}) is much more involved
because of the presence of different eigenvalues with absolute values
bigger than one. We have the following general estimate 
\begin{thm}
\label{thmart}
Let $A\in SL(d,\IZ)$ be ergodic.
There exist constants $C_{1}$ and $C_{2}$ such that for any $0<\del<1$ and sufficiently 
large $n$
\bea
\label{ulestimateA}
C_{1} e^{(1-\del)2\alpha\tilde{h}(A) n} \leq
\min_{0 \not=\bk \in \IZ^{d}} \sum_{l=1}^{n}|A^{l}\bk|^{2\alpha} \leq
C_{2}ne^{2\alpha \hat{h}(A)n}
\eea
where as before $\hat{h}(A)$ denotes minimal dimensionally averaged entropy of $A$
and $\tilde{h}(A)$ denotes a constant satisfying $0<\tilde{h}(A)\leq \hat{h}(A)$,
with equality achieved for all diagonalizable matrices $A$.
\end{thm}
The question whether the equality $\tilde{h}(A)=\hat{h}(A)$ holds for all ergodic 
matrices remains open.

The proof of the theorem relies on nontrivial use of three number-theoretical results stated below.

\textbf{I. Minkowski's Theorem on linear forms}

{\em Let $L_{1},...,L_{d}$ be linearly independent linear forms on $\IR^{d}$ which are real
or occur in conjugate complex pairs. Suppose $a_{1},a_{2},...,a_{d}$ are real positive 
numbers satisfying $a_{1}a_{2}...a_{d}=1$ and $a_{i}=a_{j}$, whenever $L_{i}=\bar{L}_{j}$. 
Then there exists a nonzero integer vector
$\bk \in \IZ^{d}$ such that for every $j=1,...,d$,
\bea
\label{minko}
|L_{j}\bk| \leq Da_{j},
\eea
where $D=|\det[L_{1},...,L_{d}]|^{1/d}$.
}

Minkowski's Theorem on linear forms will be used to obtain a sharp upper bound
on the asymptotic solution of the arithmetic minimization problem.
The proof of the above theorem and its generalization to arbitrary lattices can be
found in \cite{Newman} (Chap. VI). 

\textbf{II. Schmidt's Subspace Theorem}

{\em Let $L_{1},...,L_{d}$ be  linearly independent linear forms on $\IR^{d}$ with real
or complex algebraic coefficients. Given $\del>0$, there are finitely many proper
rational subspaces of $\IR^{d}$ such that every nonzero integer vector $\bk$ with 
\bea
\label{Schmi}
\prod_{j=1}^{d}|L_{j}\bk|<|\bk|^{-\del}
\eea   
lies in one of these subspaces.
}

Schmidt's Subspace Theorem will be used in conjunction with
Van der Waerden's Theorem on arithmetic progressions (see below)
to obtain a sharp lower bound for the asymptotic solution
of the arithmetic minimization problem.
The proof of
Schmidt's Subspace Theorem can be found in \cite{Schmidt} (Theorem 1F, p. 153).

\begin{defin}
\label{exept}
For a given set of linear forms and for fixed $\del>0$, 
the smallest collection of proper rational subspaces of $\IR^{d}$ which
contain all nonzero integer vectors satisfying (\ref{Schmi}), is called the exceptional 
set and denoted by $E_{\del}$. 
\end{defin}

A main difficulty to be resolved in using Schmidt's Subspace Theorem
is to show that the minimizer of either the original problem (\ref{artmin0})
or an equivalent problem does not lie in the respective exceptional set which
is in general unknown.
We will pursue the latter route by using
Van der Waerden's Theorem on arithmetic progressions to
show that one can always construct an equivalent minimization problem
whose minimizer is guaranteed to lie outside the corresponding exceptional
set.
To this end
we note that Schmidt's Subspace Theorem is true when the standard lattice $\IZ^{d}$ 
is replaced by any other rational lattice, that is any 
lattice of the form $\Lambda=Q(\IZ^{d})$ where $Q\in GL(d,\IQ)$. 
Schmidt's subspace theorem  can be generalized  to this situation
by considering the set of new forms 
$\tilde{L}_{j}=L_{j}Q$. The fact that $Q \in GL(d,\IQ)$ 
implies immediately that $\tilde{L}_{j}$ are still linearly independent forms on 
$\IR^{d}$ with real or complex algebraic coefficients. 

\textbf{III. Van der Waerden's Theorem on arithmetic progressions}

{\em Let $k$ and $d$ be two arbitrary natural numbers. Then there exists a
natural number $n_*(k,d)$ such that, if an arbitrary segment of length $n\geq n_*$
of the sequence of natural numbers is divided in any manner into $k$
(finite) subsequences, then an arithmetic progression
of length $d$ appears in at least one of these subsequences. 
}

The original proof was published in \cite{vdW};
Lukomskaya's simplification can be found in \cite{Lukomskaya}.

Before presenting the proof of our main results
we state a number of technical facts
concerning the structure of toral automorphisms.

\subsection{Algebraic structure of toral automorphisms}
\label{sAA}
In this section we denote by $GL(d,\IQ)$ the group of nonsingular $d\times d $ matrices with
rational entries
or the group of linear operators on Euclidean space $\IR^{d}$, which are represented in standard 
basis by such matrices. We generally use the same symbol to denote both operator and its matrix.  

In the sequel a vector $x \in \IR^{d}$ will be called an integer (or integral) vector if all 
its components are integers, and similarly a rational, an algebraic vector if all its 
components are rational or respectively algebraic numbers.
The term {\em rational subspace of $\IR^{d}$} will then refer to a linear subspace of $\IR^{d}$
spanned by rational vectors (cf. \cite{Schmidt} p. 113).  

\begin{defin}
$A \in GL(d,\IQ)$ is called irreducible (over $\IQ$) if its characteristic polynomial 
is irreducible in $\IQ[x]$.
\end{defin}

\begin{lem}
\label{irred}
The following statements about a matrix $A \in GL(d,\IQ)$ are equivalent.

a) $A$ is irreducible.

b) $A$ does not possess any proper rational $A$-invariant subspaces of $\IR^{d}$.

c) No rational proper subspace of $\IR^{d}$ is contained in any proper 
   $A$-invariant subspace of $\IR^{d}$.

d) For any nonzero $\bq \in \IQ^{d}$ and any arithmetic progression of integer 
   numbers $n_{1},...,n_{d}$, the set 

   $\{A^{n_{1}}\bq,A^{n_{2}}\bq,...,A^{n_{d}}\bq\}$ forms a basis of $\IR^{d}$.

e) $A^\dagger$ is irreducible. 

f) No nonzero $\bq \in \IQ^{d}$ is orthogonal to any proper $A$-invariant 
   subspace of $\IR^{d}$.

g) No proper $A$-invariant subspace of $\IR^{d}$ is contained in any proper
   rational subspace of $\IR^{d}$.

\end{lem}

\begin{defin}
We say that operator $A \in GL(d,\IQ)$ is completely decomposable over $\IQ$ if there 
exists a rational basis of $\IR^{d}$ in which $A$ admits the following block diagonal form  
\bea
\label{blockdiag}
\begin{bmatrix}
A_{1} & 0     & ...& 0    \\
0     & A_{2} & ...& 0    \\
...   & ...   & ...& ...  \\
0     & 0     & ...& A_{r}   
\end{bmatrix}, 
\eea
where for each $j=1,...,r \leq d$, $A_{j} \in GL(d_{j},\IQ)$ is irreducible and 
$\sum_{j=1}^r d_{j}=d$.
\end{defin}

In general, any matrix $A \in GL(d,\IQ)$ admits a rational block diagonal 
representation $[A_{j}]_{j=1,...,r}$. The smallest rational blocks to 
which $A$ can be decomposed are called elementary divisor blocks.
The characteristic polynomial corresponding to any elementary divisor block is of the form 
$p^{m}$, where $p$ is an irreducible (over $\IQ$) polynomial
 (see, e.g., \cite{Dummit}).
Although elementary divisor blocks cannot be decomposed over $\IQ$ into smaller invariant blocks, 
some elementary divisor blocks may not be irreducible. This happens
iff $m>1$ iff
$A$ is not completely decomposable over $\IQ$. One has the following
elementary fact (see Appendix \ref{aP} for a proof).  
\begin{prop}
\label{red}
 $A \in GL(d,\IQ)$ is completely decomposable over $\IQ$ iff $A$ is diagonalizable.
\end{prop}

However, even if $A\in GL(d,\IQ)$ is not completely decomposable, each elementary divisor block 
of $A$ can be uniquely represented (in a rational basis) in the following block upper triangular form
\bea
\label{but}
\begin{bmatrix}
B & C  \\
0 & D    
\end{bmatrix}, 
\eea
where $B$ is the unique rational irreducible sub-block associated with $A$-invariant
rational subspace of that elementary divisor and $C$, $D$ denote some rational matrices.

\begin{prop}
\label{distinct}
All the eigenvalues of an irreducible matrix $A\in GL(d,\IQ)$ are
distinct (complex) algebraic numbers. In particular all irreducible matrices are diagonalizable.
\end{prop}
The proofs of the above  propositions can be found in Appendix B.

\commentout{
\begin{defin}
We say that matrix $A\in GL(d,\IQ)$ is indecomposable if it does not possesses any proper 
elementary divisor blocks.
\end{defin}

It is easy to construct a nonergodic toral automorphism which has positive entropy. 
One can simply consider a toral automorphism composed of two invariant blocks:
one nonergodic and the other one ergodic. Obviously the resulting automorphisms is not 
indecomposable. The question whether there exists indecomposable nonergodic toral 
automorphism of positive entropy has a negative answer, as the following
proposition shows. 

\begin{prop}
\label{underg}
Let $F$ be an indecomposable toral automorphism. Then $F$ has positive entropy 
if and only if $F$ is ergodic.
\end{prop}

For the proofs of the above propositions we refer the reader to Appendix B.

\begin{cor}
 Let $F$ be an indecomposable and nonergodic toral automorphism. Then all eigenvalues
of $F$ are roots of unity. In particular its entropy is zero.
\end{cor}
}

Finally we note that since the leading coefficient and constant term of a characteristic 
polynomial of any toral automorphism are equal to 1, the only possible rational eigenvalues 
of such map are $\pm 1$ or $\pm i$. The latter fact implies that ergodic toral automorphisms 
do not possesses rational eigenvalues. Thus we have the following

\begin{cor}
\label{23dim}
 Let $F$ be an ergodic, two or three dimensional toral automorphism. Then $F$ is 
 irreducible (and hence diagonalizable). 
\end{cor}

\subsection{Proof of Theorem \ref{thmart}}
\label{sAP}
This section is entirely devoted to the proof of Theorem \ref{thmart}.
 
Let $[A_{j}]_{j=1,...,r}$ be a rational block-diagonal decomposition 
of $A$ into elementary divisor blocks. Since $A\in SL(d,\IZ)$, there exist 
a transition matrix $Q\in SL(d,\IQ)$ such that for every $l\in\IZ$, 
\be
A^{l}=Q^{-1}([A_{j}])^{l}Q
\ee
and moreover each elementary divisor block $[A]_{j}$ is represented in its 
block upper triangular form (\ref{but}).   

The matrix $Q$ defines a new lattice $\Lambda=Q(\IZ^{d})$ and acts bijectively between 
this lattice and the standard lattice $\IZ^{d}$. Hence   
\bean 
\min_{0 \not=\bk \in \IZ^{d}} \sum_{l=1}^{n}|A^{l}\bk|^{2\alpha}=
\min_{0 \not=\bk \in \IZ^{d}} \sum_{l=1}^{n}|Q^{-1}([A_{j}])^{l}Q\bk|^{2\alpha}=
\min_{0 \not=\bq \in \Lambda} \sum_{l=1}^{n}|Q^{-1}([A_{j}])^{l}\bq|^{2\alpha}.
\eean      
Moreover
\bean
\|Q\|^{-2\alpha}|([A_{j}])^{l}\bq|^{2\alpha}\leq
|Q^{-1}([A_{j}])^{l}\bq|^{2\alpha}\leq
\|Q^{-1}\|^{2\alpha}|([A_{j}])^{l}\bq|^{2\alpha}, \quad\forall l, j, \alpha.
\eean
Now we decompose $\Lambda$ into the direct sum of lower dimensional sublattices $\Lambda_{j}$
corresponding to invariant blocks $[A_{j}]$. So that
\be
\label{minimal}
\min_{0 \not=\bq \in \Lambda} \sum_{l=1}^{n}|([A_{j}])^{l}\bq|^{2\alpha}=
\min_{j\in \{1,...,r\}}\min_{0 \not=\bq \in \Lambda_{j}} \sum_{l=1}^{n}|(A_{j})^{l}\bq|^{2\alpha}.
\ee
Thus, without loss of generality, we may  specialize to the case that $A$ is 
already indecomposable over $\IQ$ i.e.  $A$ 
does not possesses any proper elementary divisor blocks.
To simplify the notation we will work with the standard lattice 
$\Lambda=\IZ^{d}$. According to  the
remarks following the statements of Minkowski's and Schmidt's 
Theorems  the proof can be easily adapted for
 any rational lattice $\Lambda=Q(\IZ^{d})$. 

Since the technique of the proof differs depending on diagonalizability of $A$
we consider two cases:
\subsubsection{Diagonalizable case}

Here we concentrate on the case when $A$ is diagonalizable and hence due 
to its in-decomposability irreducible (cf. Proposition \ref{red}).

We denote by $\lam_{j}$ ($j=1,...,d$) the eigenvalues of $A$.
Following Proposition \ref{distinct} we note that $\lam_{j}$ are distinct (possibly complex) algebraic 
numbers and hence there exists a basis (of $\IC^{d}$) $\{\bv_{j}\}_{j=1,...,d}$ composed of 
normalized algebraic eigenvectors corresponding to eigenvalues $\lam_{j}$.

We denote by $[P_{j}]_{j=1}^d$ the projections on $[\bv_{j}]$, and by $[L_{j}]$ the 
corresponding linear forms. It is easy to check that $[L_{j}]$ are given, in the Riesz identification, by 
the eigenvectors $[\bu_j]$ of the matrix $A^\dagger$ which are co-orthogonal to
$[\bv_j]$, i.e., $\la \bu_i,\bv_j\ra=0$ for $i\not=j$.
$[\bu_j]$ and $[\bv_j]$ are real or
occur in complex conjugate pairs. We have
\bean
 \bx=\sum_{j=1}^{d}P_{j}\bx=\sum_{j=1}^{d}(L_{j}\bx)\bv_{j}
=\sum_{j=1}^d\la \bx,\bu_j\ra \bv_j, \quad\forall\bx\in \IR^d.
\eean  

The equivalence between any two norms in a finite dimensional vector space, implies the existence 
of absolute constants $C_{1},C_{2}$ such that
\bean
 C_{1} \sum_{j=1}^{d}|P_{j}\bx|^{2} \leq |\bx|^{2} \leq C_{2} \sum_{j=1}^{d}|P_{j}\bx|^{2}.
\eean
Using the above inequalities, the monotonicity of a map $\bx \mapsto \bx^{\alpha}$ and 
an obvious inequality $(a+b)^{\alpha} \leq a^{\alpha}+ b^{\alpha}$, which holds for all 
positive $a,b$ and $\alpha \in (0,1]$ one obtains the following estimates
\bean
\sum_{l=1}^{n}|A^{l}\bk|^{2\alpha} &\leq&
\sum_{l=1}^{n}\left(C_{2}\sum_{j=1}^{d}|P_{j}A^{l}\bk|^{2}\right)^{\alpha}=
C_{2}^{\alpha}\sum_{l=1}^{n}\left(\sum_{j=1}^{d}|\lam_{j}|^{2l}|P_{j}\bk|^{2}\right)^{\alpha}\\
&\leq& C_{2}^{\alpha}\sum_{l=1}^{n}\sum_{j=1}^{d}|\lam_{j}|^{2\alpha l}|P_{j}\bk|^{2\alpha}=
C_{2}^{\alpha}\sum_{j=1}^{d}\left(\sum_{l=1}^{n}|\lam_{j}|^{2\alpha l}\right)|P_{j}\bk|^{2\alpha}
\eean
and on the other hand
\bean
\sum_{l=1}^{n}|A^{l}\bk|^{2\alpha} &\geq&
\left(\sum_{l=1}^{n}|A^{l}\bk|^{2}\right)^{\alpha} \geq
\left(\sum_{l=1}^{n} C_{1}\sum_{j=1}^{d}|P_{j}A^{l}\bk|^{2}\right)^{\alpha}\\
&=&C_{1}^{\alpha}\left(\sum_{l=1}^{n}\sum_{j=1}^{d}|\lam_{j}|^{2l}|P_{j}\bk|^{2}\right)^{\alpha}=
C_{1}^{\alpha}\left(\sum_{j=1}^{d}\left(\sum_{l=1}^{n}|\lam_{j}|^{2l}\right)|P_{j}\bk|^{2}\right)^{\alpha}.
\eean
Now we introduce some notation 
\bea
\label{lamhat}
\hat{\lam}_{j}&:=&max\{1,|\lam_{j}|\},\\
\label{hatlamgeo}
\hat{\lam}_{geo}&:=&\left(\prod_{j=1}^{d}\hat{\lam}_{j} \right)^{1/d}. 
\eea

One can easily observe that there exists a constant $C$ such that
\bean
 C\hat{\lam}_{j}^{2\alpha n}\leq \sum_{l=1}^{n}|\lam_{j}|^{2\alpha l}
 \leq n\hat{\lam}_{j}^{2\alpha n}.
\eean
In the sequel we do not distinguish between particular values of constants appearing 
in computations. The symbols $C_{1},C_{2},..$ are used to denote any generic
constants independent of $n$.

The normalization condition $|\bv_j|=1$ implies the following relation
\bea
\label{normalization}
|P_{j}\bx|=|L_{j}\bx|. 
\eea
Combining 
the above estimates one gets the following general bounds
\bea
\label{ulestimate}
C_{1}\left(\sum_{j=1}^{d} \hat{\lam}_{j}^{2 n}|L_{j}\bk|^{2}\right)^{\alpha} \leq
\sum_{l=1}^{n}|A^{l}\bk|^{2\alpha} \leq
C_{2}n\sum_{j=1}^{d} \hat{\lam}_{j}^{2\alpha n}|L_{j}\bk|^{2\alpha}.
\eea
Therefore in order to estimate (\ref{artmin0}) it suffices, essentially, to estimate
\be
\label{artmin2}
\min_{0 \not=\bk \in \IZ^{d}} \sum_{j=1}^{d}\hat{\lam}_{j}^{2\alpha n}|L_{j}\bk|^{2\alpha}.
\ee
We denote by ${\bz}_{n}$ the sequence of minimizers i.e. nonzero integral vectors
solving (\ref{artmin2}).

\medskip
\textbf{Upper bound.}
For the upper bound we assign to the set of linear forms $L_{j}$ the set $\mathcal{A}$ 
composed of all real vectors $\ba=(a_{1},...,a_{d})$ satisfying the conditions
$a_{j}>0$, for $j=1,...,d$ and $a_{i}=a_{j}$ whenever $L_{i}=\bar{L}_{j}$ and 
\be
\label{constr}
\prod_{j=1}^{d}a_{j}=1.
\ee
>From Minkowski's theorem on linear forms, we know that
for any $\ba\in\mathcal{A}$, there exists nonzero integral vector $\bk_{\ba}$ satisfying 
$|L_{j}\bk_{\ba}| \leq Da_{j}$, $j=1,...,d$, where $D=|\det[L_{1},...,L_{d}]|^{1/d}$.
 
Thus
\bea
\label{up1}
\sum_{j=1}^{d}\hat{\lam}_{j}^{2\alpha n}|L_{j}\bk_{\ba}|^{2\alpha} \leq 
D\sum_{j=1}^{d}\hat{\lam}_{j}^{2\alpha n}a_{j}^{2\alpha}.
\eea
The minimizing property of $\bz_{n}$ implies that for 
any $\ba \in \mathcal{A}$,
\bea
\label{up3}
\sum_{j=1}^{d}\hat{\lam}_{j}^{2\alpha n}|L_{j}{\bz}_{n}|^{2\alpha} \leq 
\sum_{j=1}^{d}\hat{\lam}_{j}^{2\alpha n}|L_{j}\bk_{\ba}|^{2\alpha}.
\eea
Thus combining (\ref{up1}) and (\ref{up3}), and applying the Lagrange multipliers 
minimization with the constraint (\ref{constr}) (and using the fact
that $\hat{\lam}_{i}=\hat{\lam}_{j}$ whenever $L_{i}=\bar{L}_{j}$), we get
\bea
\label{ubound1}
\sum_{j=1}^{d}\hat{\lam}_{j}^{2\alpha n}|L_{j}{\bz}_{n}|^{2\alpha} \leq 
D\min_{\ba\in\mathcal{A}}\sum_{j=1}^{d}\hat{\lam}_{j}^{2\alpha n}a_{j}^{2\alpha}=
dD\left(\prod_{j=1}^{d}\hat{\lam}_{j}^{2\alpha n}\right)^{1/d}=
dD\hat{\lam}_{geo}^{2 \alpha n}.
\eea 
Thus the
following upper bound holds
\be
\label{ubound}
 \min_{0 \not=\bk \in \IZ^{d}} \sum_{l=1}^{n}|A^{l}\bk|^{2\alpha}\leq
 C_{2} n\hat{\lam}_{geo}^{2 \alpha n}.
\ee

\medskip
\textbf{Lower bound.}
Let $m$ denote an arbitrary natural number. Using the fact that $A$ acts
bijectively on $\IZ^{d}$ we can restate the minimization problem (\ref{artmin2}) 
in the following form 
\bea
\label{artmin3}
\min_{0 \not=\bk \in \IZ^{d}}  \sum_{j=1}^{d}\hat{\lam}_{j}^{2\alpha n}|L_{j}\bk|^{2\alpha} &=&
\min_{0 \not=\bk \in \IZ^{d}}  
\sum_{j=1}^{d}\hat{\lam}_{j}^{2\alpha n}|L_{j}A^{-m}A^{m}\bk|^{2\alpha} \\
&=&\min_{0 \not=\bk \in \IZ^{d}}  
\sum_{j=1}^{d}\hat{\lam}_{j}^{2\alpha n}|\lam_{j}|^{-2\alpha m}|L_{j}A^{m}\bk|^{2\alpha}
\eea
That is
\bea
\label{eqref3}
\sum_{j=1}^{d}\hat{\lam}_{j}^{2\alpha n}|L_{j}{\bz}_{n}|^{2\alpha}=
\sum_{j=1}^{d}\hat{\lam}_{j}^{2\alpha n}|\lam_{j}|^{-2\alpha m}|L_{j}A^{m}\bz_{n}|^{2\alpha}.
\eea
We choose arbitrary $\del>0$ and consider the exceptional set $E_{\del}$ 
(see Definition \ref{exept}) associated with the system of linear forms $[L_{j}]$. 
Since $[L_{j}]$ correspond to the eigen-pairs $[\bar{\lam}_j, \bu_j]$  of $A^\dagger$
they are linearly
independent linear forms with (real or complex) algebraic coefficients. Thus the subspace theorem 
asserts that $E_{\del}$ is a finite collection of proper rational subspaces of $\IR^{d}$.   
We denote by $k_{\del}$ the number of subspaces forming $E_{\del}$. 

Now we want to show that for all sufficiently large $n$ there exist an integer
$m\leq n$ such that $A^{m}\bz_{n}$ 
does not lie in any element of $E_{\del}$. 
To this end we assume to the contrary that all $A^{m}\bz_{n}$ lie in the subspaces forming $E_{\del}$ and we 
divide the sequence of natural numbers $1,...,n$ into $k_{\del}$ classes in such 
a way that two numbers $m_{1}$ and $m_{2}$ are in the same class if $A^{m_{1}}\bz_{n}$ and 
$A^{m_{2}}\bz_{n}$ lie in the same element of $E_{\del}$.
Now let $n_*(k_{\del},d)$ be the number given in the van der Waerden theorem and let
$n\geq n_*$. 
Then there exists an arithmetic progression $m_{1},...,m_{d}$ in one of these subsequences.
By Lemma \ref{irred} d) the set of vectors
$\{A^{m_{1}}\bz_{n},A^{m_{2}}\bz_{n},...,A^{m_{d}}\bz_{n}\}$ 
forms a basis of the whole space $\IR^{d}$,
which contradicts the fact that they lie in one fixed rational proper subspace.
Hence for any $\del>0$ and $n\geq n_*$ there exists $m_*\leq n$ such that 
$A^{m_*}{\bz}_{n}$ does not lie in any element of  $E_{\del}$.

Now, introducing the notation
\be
\label{minrel1}
\hat{\bz}_{n}=A^{m_*}{\bz}_{n}
\ee
one concludes from (\ref{eqref3}) that for any $\del >0$ and all $n \geq n_*$
the following equality and estimate hold
\bea
\label{minrel2}
\sum_{j=1}^{d}\hat{\lam}_{j}^{2\alpha n}|L_{j}{\bz}_{n}|^{2\alpha}&=&
\sum_{j=1}^{d}\hat{\lam}_{j}^{2\alpha n}|\lam_j|^{-2\alpha m_{*}}|L_{j}\hat{\bz}_{n}|^{2\alpha}\\
\label{constr1}
\prod_{j=1}^{d}|L_{j}\hat{\bz}_{n}| &\geq& \frac{1}{|\hat{\bz}_{n}|^{\del}}.
\eea
Inequality (\ref{constr1}) may be rewritten as 
\be
\label{f}
  \prod_{j=1}^{d} |L_{j}\hat{\bz}_{n}|=\frac{1}{f(|\hat{\bz}_{n}|)^{\del}}
\ee
with some $f:\IR^+\to\IR^+$ such that $f(r)\leq r, \forall r>0$.

Using (\ref{minrel1}) and (\ref{ubound1}) we obtain the existence of 
a constant $\lam>1$ such that
\bea
\label{minbound}
f(|\hat{\bz}_{n}|) \leq 
|\hat{\bz}_{n}|=
|A^{m_*}{\bz}_{n}| \leq 
\hat{\lam}_{max}^{m_*}|{\bz}_{n}| \leq 
\hat{\lam}_{max}^{n}\sum_{j=1}^{d}\hat{\lam}_{j}^{n}|L_{j}{\bz}_{n}| \leq
dD(\hat{\lam}_{max}\hat{\lam}_{geo})^{n} \leq
\lam^{n}.
\eea
Note that $\prod_{j}\lam_j=1$. So, by (\ref{f}) the quantities 
$B_{j,n}=\left(|\lam_{j}|^{-m_*}f(|\hat{\bz}_{n}|)^{\del/d}|L_{j}\hat{\bz}_{n}|\right)^{2\alpha}, j=1,...,d$ 
satisfy the constraint
\bea
\label{newc}
\prod_{j=1}^{d} B_{j,n}=1,\quad\forall n>n_*.
\eea
Thus applying (\ref{minbound}) and the Lagrange multipliers minimization with the
constraint  (\ref{newc}) one gets
\bean
\sum_{j=1}^{d}\hat{\lam}_{j}^{2\alpha n}|\lam_{j}|^{-2\alpha m_{*}}|L_{j}\hat{\bz}_{n}|^{2\alpha}
= f(|\hat{\bz}_{n}|)^{-2\alpha\del/d}\sum_{j=1}^{d}\hat{\lam}_{j}^{2\alpha n}B_{j,n}\geq 
  \lam^{-2 \alpha n\del/d}\hat{\lam}_{geo}^{2\alpha n}=:  
  \hat{\lam}_{geo}^{2\alpha n(1-\hat{\del})}.
\eean
This and equality (\ref{minrel2}) yields the following lower bound for (\ref{artmin2})
\bea
\label{dlbound}
\min_{0 \not=\bk \in \IZ^{d}} \sum_{j=1}^{d}\hat{\lam}_{j}^{2\alpha n}|L_{j}\bk|^{2\alpha} \geq 
\hat{\lam}_{geo}^{2 \alpha n(1-\hat{\del})}.
\eea

\subsubsection{Non-diagonalizable case}

We move on to the general case where $A$ is not irreducible (but,
as assumed at the beginning of the proof, indecomposable over $\IQ$).
We denote by $B$ the invariant irreducible sub-block of $A$ given by its block upper triangular
decomposition (\ref{but}) and by $S$ the rational invariant subspace associated with this block.
We note that $B$ as an irreducible matrix is diagonalizable.

\medskip
\textbf{Upper bound.}
Note that
\bea
 \min_{0 \not=\bk \in \IZ^{d}} \sum_{l=1}^{n}|A^{l}\bk|^{2\alpha} \leq 
 \min_{0 \not=\bk \in S \cap \IZ^{d}} \sum_{l=1}^{n}|B^{l}\bk|^{2\alpha}.
\eea
The corresponding upper bound (\ref{ubound}) for $B$ is then also an upper 
bound for the whole matrix $A$. We note that geometric average of $\hat{\lam}_j$ over $S$ is
equal to the geometric average of all $\hat{\lam}_j$ associated with matrix $A$ 
(i.e. over the whole space $\IR^{d}$).

\medskip
\textbf{Lower bound.}
According to our assumption $A$ is indecomposable and thus the characteristic polynomial
of $A$ is of the form $p^m$ for some irreducible $p$. All 
Jordan blocks of $A$ have the same size $m$ and different Jordan blocks correspond to distinct eigenvalues.
Denote by $b$ the number of the Jordan blocks in $A$ and 
by $\lam_{j}$, where $j=1,..,l$ all these distinct eigenvalues . Since each $\lam_{j}$ has algebraic multiplicity 
$m$, we get $d=mb.$
Let $\{\bv_{j,h}\}_{j=1,...,b;h=0,...,m-1}$ be a basis (of $\IC^{d}$) 
in which $A$ admits the Jordan canonical form. As usually $L_{j,h}$ will denote the corresponding linear forms.
Each  $\bv_{j,h}$ can be regarded as a generalized eigenvector of $A$ associated with 
an eigenvalue $\lam_{j}$. We assume that these generalized eigenvectors are ordered according to
their degree i.e. $\bv_{j,h}$ satisfies the equation $(A-\lam_{j}I)^{1+h}\bv_{j,h}=0$.   
Reordering the eigenvalues, if necessary, we can also assume that $\lam_{1}$ has the largest modulus among
all eigenvalues of $A$ and hence $\hat{\lam}_{1}=|\lam_{1}|$.
Let ${\bz}_{n}$ be the sequence of minimizers solving (\ref{artmin0}).
We first note that for each $n$ there exists 
$0 \leq h \leq m-1$ such that $L_{1,h}{\bz}_{n} \not =0$. 
Indeed, otherwise for all $h=0,...,m-1$, $L_{1,h}{\bz}_{n} =0$ and
consequently for any $n$ and $h$ $L_{1,h}A^{n}{\bz}_{n}=0$. The latter implies that the set of consecutive
iterations $\{{\bz}_{n},A^{1}{\bz}_{n},A^{2}{\bz}_{n},...\}$ spans a proper rational $A$-invariant 
subspace of $\IR^{d}$ which does not have any intersection with the subspace spanned by the generalized eigenvectors 
of $A$ associated with eigenvalue $\lam_{1}$. This clearly contradicts the irreducibility of $p$. 
Now, for given $n$ we denote by $h(n)$ the biggest index $h$ for which the condition $L_{1,h}{\bz}_{n} \not =0$
holds.

We have the following estimate
\bea
\label{est1}
 &&\hat{\lam}_{1}^{2\alpha n}|L_{1,h(n)}{\bz}_{n}|^{2\alpha} \leq
 \left(\sum_{j=1}^{b}\sum_{h=0}^{m-1}
 \left|\sum_{i=0}^{m-1-h}\lam_{j}^{n-i}\binom{n}{i}L_{j,h+i}{\bz}_{n}\right|^{2}\right)^{\alpha}\\
 &\leq&C_{1}|A^{n}{\bz}_{n}|^{2\alpha} \leq
 C_{1}\sum_{l=1}^{n}|A^{l}{\bz}_{n}|^{2\alpha} \leq
 C_{2}n\hat{\lam}_{geo}^{2\alpha n},
\eea
where the last inequality follows from previously established upper bound.

>From the Diophantine approximation and the assumption that $|L_{1,h(n)}{\bz}_{n}|\not=0$, 
there exists $\beta>0$ such that (see \cite{Schmidt} p. 164) 
\bea
\label{est2}
|L_{1,h(n)}{\bz}_{n}| \geq \frac{1}{|{\bz}_{n}|^{\beta}}.
\eea 
Thus combining (\ref{est1}) with (\ref{est2}) one gets
\bean
 \hat{\lam}_{1}^{2\alpha n}|{\bz}_{n}|^{-2\alpha\beta} \leq
 \hat{\lam}_{1}^{2\alpha n}|L_{1,h(n)}{\bz}_{n}|^{2\alpha} \leq 
 C_{2}n\hat{\lam}_{geo}^{2\alpha n}.
\eean

After rearrangements one obtains the following 
lower bound estimate for (\ref{artmin0}) 
\bea
\label{nlbound}
\frac{C}{n}\hat{\lam}^{2\alpha n} \leq 
C|{\bz}_{n}|^{2\alpha} \leq 
\sum_{l=1}^{n}|A^{l}{\bz}_{n}|^{2\alpha}=
\min_{0 \not=\bk \in \IZ^{d}} \sum_{l=1}^{n}|A^{l}\bk|^{2\alpha},
\eea
where
\be
 \hat{\lam}=\left(\frac{\hat{\lam}_{1}}{\hat{\lam}_{geo}}\right)^{1/\beta}.
\ee
We note that ergodicity of $A$ implies $\hat{\lam}_{1}>\hat{\lam}_{geo}>1$ 
(see (\ref{hatlamgeo}), (\ref{entropy}) and Proposition \ref{zeroentropy}) which ensures 
non-triviality of this lower bound.

Now in order to finish the proof is suffices to combine the estimates (\ref{ubound}),
(\ref{dlbound}) and (\ref{nlbound}), and note that 
\bean
\hat{\lam}_{geo}^{2\alpha n}=e^{2\alpha\frac{h(A)}{d}n}=e^{2\alpha\hat{h}(A)n}
\eean
which yields (\ref{ulestimateA}). $\qquad \blacksquare$

\subsection{Proofs of Theorem \ref{thm2} ii), iii) and Theorem \ref{thm3}}
\label{sAT}
In this section we apply Theorem \ref{thmart} to prove main theorems of Section \ref{sD}.

In order to determine the dissipation time of $T_{\vep,\alpha}$ one has to 
determine the asymptotics of $\|T_{\vep,\alpha}^{n}\|$ when $n$ goes to 
infinity. 
According to formulas (\ref{Tnorm3}) and (\ref{artminA}) this problem reduces 
to problem 
(\ref{artmin0}) solved in previous sections.

Thus in view of Theorem \ref{thmart})
there exist constants $C_{1}$ and $C_{2}$ such that for any $\del,\del'>0$ and sufficiently 
large $n$
\bean
\label{ulestimate1}
C_{1} e^{(1-\del)2\alpha\tilde{h}(A) n} \leq
\min_{0 \not=\bk \in \IZ^{d}} \sum_{l=1}^{n}|A^{l}\bk|^{2\alpha} \leq
C_{2}ne^{2\alpha \hat{h}(A)n}\leq
C_{2}e^{(1+\del')2\alpha \hat{h}(A)n}
\eean
Using formula (\ref{Tnorm3})
\bean
e^{-\vep  C_{2}e^{(1+\del')2\alpha \hat{h}(A)n}}  
\leq \|T_{\vep,\alpha}^{n}\| \leq 
e^{-\vep  C_{1}e^{(1-\del)2\alpha\tilde{h}(A)n}}.
\eean
Now when $n=n_{diss}$, we have
\bean
C_{1} e^{(1-\del)2\alpha\tilde{h}(A)n_{diss}}
\leq \frac{1}{\vep} \leq
C_{2} e^{(1+\del')2\alpha \hat{h}(A)n_{diss}}  
\eean
and 
\bean
\frac{1}{(1+\del')2\alpha\hat{h}(A)}\bigg(\ln(1/\vep) -\ln C_{2}\bigg)   
\leq n_{diss} \leq
\frac{1}{(1-\del)2\alpha \tilde{h}(A)}\bigg(\ln(1/\vep)-\ln C_{1}\bigg),
\eean
which proves part ii) of Theorem \ref{thm2} i.e. the logarithmic growth of 
dissipation time as a function of $1/\vep$. 

Moreover, using the definition of dissipation rate constant
\bean
  R_{diss}=\lim_{\vep \rightarrow 0} \frac{n_{diss}}{\ln(1/\vep)}
\eean  
we obtain
\bean
\frac{1}{(1+\del')2\alpha \hat{h}(A)} \leq 
R_{diss} \leq 
\frac{1}{(1-\del)2\alpha \tilde{h}(A)}.
\eean
Finally letting $\del \rightarrow 0$ and $\del' \rightarrow 0$  
we arrive at the following results:
\begin{itemize}
\item The general case - Theorem \ref{thm2} iii)
\bean
\frac{d}{2\alpha \hat{h}(F)}\leq R_{diss} 
\leq \frac{1}{2\alpha\tilde{h}(F) }
\eean
\item The diagonalizable case - Theorem \ref{thm3}
\bean
 R_{diss}=\frac{1}{2\alpha \hat{h}(F)}.
\eean
\end{itemize}
This completes the proof.
$\qquad \blacksquare$

\section{Degenerate Noise}
\label{degnoise}

In this section we compute the dissipation time for non-strictly contracting generalizations 
of $\alpha$-stable transition operators. 
Instead of considering standard $\alpha$-stable kernels of the form 
(\ref{kernel}) one can allow for some degree of degeneracy of noise in chosen directions by 
introducing the following family of noise kernels
\be
\label{noker}
g_{\vep,\alpha,B}(\bx):=\sum_{\bk \in \IZ^{d}}e^{-\epsilon |B\bk|^{2\alpha}}\bfe_{\bk}(\bx),
\ee
Where $B$ denotes any $d\times d$ matrix with $\det B=0$. 

We denote by $G_{\vep,\alpha,B}$ the noise operator associated with $g_{\vep,\alpha,B}$.  
The degeneracy of $B$ immediately implies that $\|G_{\vep,\alpha,B}\|=1$ and hence the
general considerations of sections 1 and 2 do not apply here.
The answer to the question whether or not the dissipation time 
is finite depends on the choice of matrix $B$. 

For simplicity we concentrate on the case when $B$ is diagonalizable.

We call the eigenvector of $B$ nondegenerate if it corresponds
to nonzero eigenvalue. 
\begin{thm}
Let $F$ be any toral automorphism and 
$T_{\vep,\alpha,B}=G_{\vep,\alpha,B}U_{F}$. 
Assume that $B$ is diagonalizable. Then
 
i) If all nondegenerate eigenvectors of $B^{*}$ lie in one proper invariant subspace of $F$ then 
dissipation does not take place i.e. $n_{diss}=\infty$.

ii) Otherwise the following statements hold.

a) $T_{\vep,\alpha}$ has simple dissipation time iff $F$ is not ergodic.

b) $T_{\vep,\alpha}$ has logarithmic dissipation time iff $F$ is ergodic. 

c) If $T_{\vep,\alpha,B}$ has logarithmic dissipation time then the dissipation rate constant satisfies the following bounds
\bean
\frac{1}{2 \alpha \hat{h}(F)} \leq R_{diss} \leq \frac{1}{2 \alpha \tilde{h}(F)},
\eean
with some constant $\tilde{h}(F)\leq\hat{h}(F)$. The equality is achieved for all diagonalizable automorphisms $F$. 
\end{thm}
\textbf{Proof.}

We continue to use the convention $A=F^{\dagger}$.
The general formula derived previously for $\|T_{\vep,\alpha}^{n}\|$ (see (\ref{Tnorm3})), will now
take the form 
\bea
\label{Tnorm3B}
\|T_{\vep,\alpha,B}^{n}\| &=& \sup_{0 \not=\bk \in \IZ^{d}} e^{-\vep \sum_{l=1}^{n}|BA^{l}\bk|^{2\alpha}} = 
e^{-\vep \underset{0 \not=\bk \in \IZ^{d}}{\inf} \sum_{l=1}^{n}|BA^{l}\bk|^{2\alpha}}. 
\eea
Thus we need to estimate
\bean
\label{Binf}
\inf_{0 \not=\bk \in \IZ^{d}} \sum_{l=1}^{n}|BA^{l}\bk|^{2\alpha}.
\eean
To this end we denote by $\mu_{j}$ ($j=1,...,d$) 
the eigenvalues of $B$
and we construct a basis (of $\IC^{d}$) $\{\bv_{j}\}_{j=1,...,d}$ composed of 
normalized eigenvectors corresponding to eigenvalues $\mu_{j}$.
We denote by ${P_{j}}_{j=1,...,d}$ the set of eigen-projections on ${\bv_{j}}$, and by ${L_{j}}$ the 
set of corresponding linear forms, given by the eigenvectors $\bu_j$ of $B^{\dagger}$, which are
of course co-orthogonal to ${\bv_j}$, i.e. $\la \bu_i,\bv_j\ra=0$ for $i\not=j$.
We have
\bean
 \bx=\sum_{j=1}^{d}P_{j}\bx=\sum_{j=1}^{d}(L_{j}\bx)\bv_{j}
=\sum_{j=1}^d\la \bx,\bu_j\ra \bv_j, \quad\forall\bx\in \IR^d.
\eean  
In subsequent computations the symbols $C_{1}$, $C_{2}$ denote some
absolute constants values of which are subject to change during calculations. 

We consider two cases.

i) All nondegenerate eigenvectors of $B^{\dagger}$ lie in one proper subspace of $F$.
We have the following estimates
\bean
|BA^{l}\bk|^{2}\geq 
C_{1} \sum_{j=1}^{d}|P_{j}BA^{l}\bk|^{2} =
C_{1} \sum_{j=1}^{d}|\mu_{j}|^{2}|P_{j}A^{l}\bk|^{2} =
C_{1} \sum_{j=1}^{d}|\mu_{j}|^{2}|\la A^{l}\bk,\bu_j\ra|^{2} =
C_{1} \sum_{j=1}^{d}|\mu_{j}|^{2}|\la \bk,F^{l}\bu_j\ra|^{2}
\eean
and 
\bean
|BA^{l}\bk|^{2}\leq 
C_{2} \sum_{j=1}^{d}|P_{j}BA^{l}\bk|^{2}=
C_{2} \sum_{j=1}^{d}|\mu_{j}|^{2}|P_{j}A^{l}\bk|^{2}=
C_{2} \sum_{j=1}^{d}|\mu_{j}|^{2}|\la A^{l}\bk,\bu_j\ra|^{2}=
C_{2} \sum_{j=1}^{d}|\mu_{j}|^{2}|\la \bk,F^{l}\bu_j\ra|^{2}
\eean
Since at least one of $\mu_{j}$ is zero and all nondegenerate vectors $\bu_j$ lie in a proper invariant 
subspace of $F$, one easily sees that for each fixed $n$ 
\bean
\inf_{0 \not=\bk \in \IZ^{d}}\sum_{l=1}^{n}|BA^{l}\bk|^{2\alpha}=
\inf_{0 \not=\bk \in \IZ^{d}}\sum_{l=1}^{n}\sum_{j=1}^{d}|\mu_{j}|^{2\alpha}|\la \bk,F^{l}\bu_j\ra|^{2\alpha}=0.
\eean
ii) In this case we have the following upper bound
\bea
\label{ubinf}
\inf_{0 \not=\bk \in \IZ^{d}} \sum_{l=1}^{n}|BA^{l}\bk|^{2\alpha}
\leq \|B\|^{2\alpha}\inf_{0 \not=\bk \in \IZ^{d}} \sum_{l=1}^{n}|A^{l}\bk|^{2\alpha}
= \|B\|^{2\alpha}\min_{0 \not=\bk \in \IZ^{d}} \sum_{l=1}^{n}|A^{l}\bk|^{2\alpha}. 
\eea
In order to provide an appropriate lower bound
we note that the set of vectors $\{F^{h}\bu_{j}\}$, where $1\leq h \leq d$ and $j$ runs through the
indices of all nondegenerate eigenvectors of $B$, spans the whole space 
(otherwise all nondegenerate $\bu_{j}$ would lie in one proper invariant subspace of $F$). 
We denote by $\{F^{h_{i}}\bu_{j_{i}}\}$ ($1\leq i \leq d$) a basis extracted from the above set.
We can define now a new norm $|\cdot|_{\bu}$ on $\IR^{d}$ by
\bean
|\bx|^{2}_{\bu}=\sum_{i=1}^{d}|\la \bx,F^{h_{i}}\bu_{j_{i}} \ra|^{2}
\eean
and compute
\bean
\sum_{l=1}^{dn}|BA^{l}\bk|^{2\alpha}
&=&\sum_{l=0}^{n-1}\sum_{h=1}^{d}|BA^{dl+h}\bk|^{2\alpha}
\geq \sum_{l=0}^{n-1}\sum_{h=1}^{d}C_{1}\sum_{j=1}^{d}|P_{j}BA^{dl+h}\bk|^{2\alpha}\\
&=&C_{1}\sum_{l=0}^{n-1}\sum_{h=1}^{d}\sum_{j=1}^{d}|\mu_{j}|^{2\alpha}|P_{j}A^{dl+h}\bk|^{2\alpha}
\geq C_{1} \sum_{l=0}^{n-1}\sum_{i=1}^{d}|\la A^{dl+h_{i}}\bk,\bu_{j_{i}} \ra|^{2\alpha}\\
&=&C_{1}\sum_{l=0}^{n-1}\sum_{i=1}^{d}|\la A^{dl}\bk,F^{h_{i}}\bu_{j_{i}} \ra|^{2\alpha}
=C_{1}\sum_{l=0}^{n-1}|A^{dl}\bk|^{2\alpha}_{\bu}
\eean
Using the equivalence between norms $|\cdot|$ and $|\cdot|_{\bu}$ and combining
(\ref{ubinf}) with the above estimate we get
\bean
C_{1}\min_{0 \not=\bk \in \IZ^{d}} \sum_{l=0}^{n-1}|A^{dl}\bk|^{2\alpha}\leq
\inf_{0 \not=\bk \in \IZ^{d}}\sum_{l=1}^{dn}|BA^{l}\bk|^{2\alpha}\leq
\|B\|^{2\alpha}\min_{0 \not=\bk \in \IZ^{d}} \sum_{l=1}^{dn}|A^{l}\bk|^{2\alpha}. 
\eean
This together with the obvious fact that $\hat{h}(A^{d})=d\hat{h}(A)$ and the general estimate
(\ref{ulestimateA}) reduces the proof back to the nondegenerate case considered in the previous section. 
$\qquad \blacksquare$

\section{Time of decay of $\vep$-coarse-grained states}
\label{coarse}

The uncertainties in the initial preparation and the final measurement of
the noiseless system give rise to non-cumulative random perturbations
to the system. Alternatively, one can coarse-grain the initial and final states of
the noiseless system by convoluting with the
$\veps$-noise kernel. 
That is, instead of the original operator $T_{\vep,\alpha}$,
we consider the operator
$\hat{T}_{\vep,\alpha}$ defined as 
\be
 \hat{T}_{\vep,\alpha}^{n}=G_{\vep,\alpha}U_{F}^{n}G_{\vep,\alpha}=G_{\vep,\alpha}U_{F^{n}}G_{\vep,\alpha}.
\ee
and compute the number of iterations required to have the $L^2$-norm of
the final state being $e^{-1}$ times that of the initial state. 
We will show that for ergodic toral automorphisms
the required number of iterations is essentially the same
asymptotically as the dissipation time computed in the previous sections.

One can represent the action of $U_{F}$ or more generally $U_{F}^{n}$
in the Fourier series 
\ben
 U_{F}^{n}\bfe_{\mf{\bk}}=\sum_{0 \not = \bk' \in \IZ^{d}}u_{\bk,\bk'}^{(n)}\bfe_{\bk'},
\een
where $u_{\bk,\bk'}^{(1)}$ coincides with $u_{\bk,\bk'}$ defined previously
(cf. (\ref{uk})) and
\[
u^{(n)}_{\bk,\bk'}=\sum_{0\neq \bk_1,...,\bk_{n-1}\in \IZ^d} u_{\bk,\bk_1}u_{\bk_1,\bk_2}...u_{\bk_{n-1},\bk'}
\]
which satisfies
\bea
\label{ul3}
\sum_{0 \not = \bk' \in \IZ^{d}}|u^{(n)}_{\bk,\bk'}|^{2}=1,\quad \forall\, n,\,\,\bk.
\eea
Then 
\bean
\hat{T}_{\vep,\alpha}^{n}\bfe_{\bk_{0}}&=&
G_{\vep,\alpha}U_F^{n}G_{\vep,\alpha}\bfe_{\bk_{0}}=
G_{\vep,\alpha}U_F^{n}e^{-\vep|\bk_{0}|^{2}} \bfe_{\bk_{0}}=
e^{-\vep|\bk_{0}|^{2}}G_{\vep,\alpha}\sum_{0 \not = \bk_{n} \in \IZ^{d}}
u^{(n)}_{\bk_{0},\bk_{n}}\bfe_{\bk_{n}}\\
&=&e^{-\vep(|\bk_{0}|^{2}+|\bk_{n}|^{2})}
\sum_{0 \not = \bk_{n} \in \IZ^{d}}u^{(n)}_{\bk_{0},\bk_{n}}\bfe_{\bk_{n}}.
\eean

Now we define
\ben
S_n(\bk_{n})=\{\bk_{0} \in \IZ^{d}\backslash \{0\} : u^{(n)}_{\bk_{0},\bk_{n}}\not = 0\}.
\een
Similar computations to these performed in Section \ref{sD}
give the following  general upper bound for $\|\hat{T}_{\vep,\alpha}^{n}\|$
\bea
\label{uppernorm1}
\|\hat{T}_{\vep,\alpha}^{n}f\|^{2} \leq
\sum_{0 \not = \bk_{n} \in \IZ^{d}}
\sum_{\bk_{0} \in S_n(\bk_{n})}
|\hat{f}(\bk_{0})|^{2}
\sum_{\bk_{0} \in S_n(\bk_{n})}
|u_{\bk_{0},\bk_{n}}^{(n)}|^{2}.
\eea
For a toral automorphism one easily sees that
\bea
u_{\bk_{0},\bk_{n}}^{(n)}=
e^{-\vep(|\bk_{0}|^{2\alpha}+|A^{n}\bk_{0}|^{2\alpha})}\del_{\bk_{n},A^{n}\bk_{0}}
\eea
and hence
\ben
\|\hat{T}_{\vep,\alpha}^{n}\| = 
e^{-\vep \min_{0 \not=\bk \in \IZ^{d}}(|\bk_{0}|^{2\alpha}+|A^{l}\bk_{0}|^{2\alpha})}. 
\een
The arithmetic minimization problem (\ref{artminA}) corresponding to
the dissipation time of $\hat{T}_{\vep,\alpha}$ now becomes
\bea 
\label{artminA1}
\min_{0 \not=\bk \in \IZ^{d}}\left(|\bk|^{2\alpha}+|A^{n}\bk|^{2\alpha}\right).
\eea
The key observation is that, by the same arguments as before,
similar estimates to these given in (\ref{ulestimate}) hold
\bea
\label{ulestimate2}
C_{1}\left(\sum_{j=1}^{d} \hat{\lam}_{j}^{2 n}|L_{j}\bk|^{2}\right)^{\alpha} \leq
|\bk|^{2\alpha}+|
A^{n}\bk|^{2\alpha} \leq
C_{2}\sum_{j=1}^{d} \hat{\lam}_{j}^{2\alpha n}|L_{j}\bk|^{2\alpha}.
\eea
The remaining computations are the same verbatim so 
the dissipation time of $T_{\vep,\alpha}$ and $\hat{T}_{\vep,\alpha}$ are 
equal asymptotically.

\section{Time Scales in Kinematic Dynamo}
\label{dynamo}
In this section we briefly discuss the connection between the dissipation time 
and some characteristic time scales associated with kinematic dynamo,
which concerns the generation of
electromagnetic fields by  
mechanical motion. 
For a general setup and discussion we refer the reader to \cite{CG} and
\cite{KY} and 
references therein.
Here we restrict ourselves only to necessary definitions.

Let $\bB\in L^{2}_{0}(\IT^{d},\IR^{d})$ denote periodic, zero mean and divergence 
free magnetic field and let $F$ be the time-1 map associated with the fluid velocity. 
We define the push-forward map 
\bean
F_{*}\bB(\bx)=dF(F^{-1}(\bx))\bB(F^{-1}(\bx)).
\eean
The noisy push-forward map $P_{\vep,\alpha}$
on $L_{0}^2(\IT^{d},\IR^{d})$ 
is then given by
\bea
\label{Tdyn}
P_{\vep,\alpha}\bB:=G_{\vep,\alpha}F_{*}\bB,
\eea   
where the convolution (the action of $G_{\vep,\alpha}$) is applied component-wise.

It is said that the kinematic dynamo action (positive dynamo effect) occurs if
the dynamo growth rate is positive i.e.
\bean
R_{dyn}=\lim_{n\rightarrow \infty}\frac{1}{n}\ln\|P_{\vep,\alpha}^{n}\|>0.
\eean
Moreover if
\bean
\lim_{\vep\rightarrow 0}\lim_{n\rightarrow \infty}\frac{1}{n}\ln\|P_{\vep,\alpha}^{n}\|>0,
\eean
then the dynamo action is said to be fast; otherwise it is slow. 
The anti-dynamo action takes place if
\bean
\lim_{n\rightarrow \infty}\frac{1}{n}\ln\|P_{\vep,\alpha}^{n}\|<0.
\eean
Now we introduce the {\em threshold time} scale as 
\bean
n_{th}=\max \{n : \|P_{\vep,\alpha}^{n}\|>e\,\, \hbox{such that}\,\,
\|P_{\vep,\alpha}^{n-1}\| \,\,\hbox{or}\,\,
\|P_{\vep,\alpha}^{n+1}\| \leq e\}.
\eean 
The threshold time $n_{th}(\vep)$ is of order $O(1)$ as $\vep\rightarrow 0$ for all
fast kinematic dynamo systems. In the case of anti-dynamo action, $n_{th}(\vep)$ captures
the longest time scale on which the generation of the magnetic field still takes place.
Finally $n_{th}(\vep)$ is not defined for systems with do not exhibit any growth of 
magnetic field throughout the evolution. 
In the case of anti-dynamo
we consider the time scale on which the generation of the
magnetic field achieves its maximal value
\bean
n_{p}=\min\{n: \|P_{\vep,\alpha}^{n}\|=\sup_{m}\|P_{\vep,\alpha}^{m}\|\}.
\eean
which is called the {\em peak time} of the anti-dynamo action.

Our next theorem establishes the relation between $n_{p}$, $n_{th}$ and $n_{diss}$ for toral automorphisms.
Thus $dF=F$ and
\bean
P_{\vep,\alpha}\bB=g_{\vep,\alpha}*F(\bB\circ F^{-1}).
\eean
\begin{thm}
Let $F$ be any toral automorphism.
Then

i) If $F$ is nonergodic and has positive entropy then for all 
$0<\vep<R_{diss}\ln\rho_{F}$ the fast dynamo action takes place
with dynamo growth rate satisfying 
\bean
R_{dyn}= \ln\rho_{F}-\vep R_{diss}^{-1}\overset{\vep\rightarrow 0}\longrightarrow \ln\rho_{F}>0,
\eean
where $\rho_{F}$ denotes the spectral radius of $F$.
The threshold time $n_{th}$ is of order $O(1)$ and 
if $F$ is diagonalizable then $n_{th}\approx [R_{dyn}^{-1}]$. 

ii) If $F$ is nonergodic and has zero entropy then anti-dynamo action occurs and 
for nondiagonalizable $F$, 
\bean
n_{p}&\sim& \frac{n_{th}}{\ln(n_{th})}\\
&\sim& n_{diss}\\
&\approx& R_{diss}\frac{1}{\vep}.
\eean
Moreover there exists a constant $0< \gamma \leq d$ such that 
$\|P_{\vep,\alpha}^{n_{p}}\|\sim (1/\vep)^{\gamma}$.
If $F$ is diagonalizable then $\|P_{\vep,\alpha}^{n}\|$ is strictly decreasing (in $n$) and, hence,
$n_{p}=0$ and $n_{th}$ is not defined.
   
iii) If $F$ is ergodic then anti-dynamo action occurs and
\bean
n_{p}\approx n_{diss}. 
\eean
In particular if $F$ is diagonalizable then
\bean
n_{p} &\approx& n_{th}-R_{diss}\ln(n_{th})\\
&\approx& R_{diss}\ln(1/\vep)\\
&=& \frac{1}{2\alpha\hat{h}(F)}\ln(1/\vep)
\eean
and 
\be
\label{anti}
\|P_{\vep,\alpha}^{n_{p}}\|\sim (1/\vep)^{\frac{\ln\rho_{F}}{2\alpha \hat{h}(F)}}.
\ee
\end{thm}
We see that even in the case of anti-dynamo action the magnetic field
can still grow to relatively large magnitude when the noise
is small (power-law in $1/\ep$).

\textbf{Proof.}
Representing the initial magnetic field $\bB=(\bb_{1},...,\bb_{d})$ in Fourier basis
\bean
\bB=\sum_{0 \not = \bk \in \IZ^{d}}\hat{\bB}(\bk)e_{\bk}
\eean
one obtains
\bean
P_{\vep,\alpha}\bB=\sum_{0 \not = \bk \in \IZ^{d}}F\hat{\bB}(\bk)
e^{-\vep |A\bk|^{2\alpha}}e_{A\bk},
\eean
where we set $A=(F^{-1})^{\dagger}$.
After $n$ iterations
\bean
P_{\vep,\alpha}^{n}\bB=\sum_{0 \not = \bk \in \IZ^{d}}F^{n}\hat{\bB}(\bk)
e^{-\vep \sum_{l=1}^{n}|A^{l}\bk|^{2\alpha}}
e_{A^{n}\bk}.
\eean
Thus
\bean
&&\|P_{\vep,\alpha}^{n}\bB\|^{2} \leq 
\sum_{0 \not = \bk \in \IZ^{d}}|F^{n}\hat{\bB}(\bk)|^{2}
e^{-2\vep \sum_{l=1}^{n}|A^{l}\bk|^{2\alpha}}
\leq \max_{0 \not=\bk \in \IZ^{d}}e^{-2\vep \sum_{l=1}^{n}|A^{l}\bk|^{2\alpha}} 
\sum_{0 \not = \bk \in \IZ^{d}}|F^{n}\hat{\bB}(\bk)|^{2}\\
&=&e^{-2\vep \min_{0 \not=\bk \in \IZ^{d}}\sum_{l=1}^{n}|A^{l}\bk|^{2\alpha}} 
|F^{n}\bB|^{2}=e^{-2\vep \sum_{l=1}^{n}|A^{l}\bk_{n}|^{2\alpha}} |F^{n}\bB|^{2}
\leq e^{-2\vep \sum_{l=1}^{n}|A^{l}\bk_{n}|^{2\alpha}} \|F^{n}\|^{2}|\bB|^{2},
\eean
where $\bk_{n}$ denotes a solution of the minimization problem 
\bean
\min_{0 \not=\bk \in \IZ^{d}}\sum_{l=1}^{n}|A^{l}\bk|^{2\alpha}.
\eean
The above calculation provides the following upper bound
\bean
\|P_{\vep,\alpha}^{n}\| &\leq& e^{-\vep \sum_{l=1}^{n}|A^{l}\bk_{n}|^{2\alpha}} \|F^{n}\|.
\eean
On the other hand let $\bv_{n}$ denote a unit vector satisfying $\|F^{n}\|=|F^{n}\bv|$. 
One immediately sees that the above upper bound for $\|P_{\vep,\alpha}^{n}\|$ is
achieved for magnetic field of the form $\bB=\bv_{n}e_{\bk_{n}}$.
Thus
\bea
\|P_{\vep,\alpha}^{n}\|= e^{-\vep \sum_{l=1}^{n}|A^{l}\bk_{n}|^{2\alpha}} \|F^{n}\|.
\eea
Now we consider the cases mentioned in the statement of the theorem

i) Nonergodic, nonzero entropy case. 

For any nonergodic map we have
\bean
\sum_{l=1}^{n}|A^{l}\bk_{n}|^{2\alpha} \approx R_{diss}^{-1} n.
\eean
This implies the following asymptotics
\bea
\label{na}
\|P_{\vep,\alpha}^{n}\|\approx e^{-\vep R_{diss}^{-1} n} \|F^{n}\|\approx 
e^{(-\vep R_{diss}^{-1} + \ln\rho_{F})n + c_{1}\ln n +c_{2}},
\eea
where $c_{1},c_{2}\geq 0$ are constants (both equal $0$ iff $F$ is diagonalizable). 
Thus for $\vep<R_{diss}{\ln\rho_{F}}$ we have
\bean
R_{dyn}= \ln\rho_{F}-\vep R_{diss}^{-1} \overset{\vep\rightarrow 0}\longrightarrow \ln\rho_{F}>0.
\eean
The threshold time is clearly of order $O(1)$ and in
diagonalizable case can be written as 
\bean
n_{th}\approx \frac{1}{\ln\rho_{F}-\vep R_{diss}^{-1}} 
\overset{\vep \rightarrow 0}\longrightarrow \frac{1}{\ln\rho_{F}}.
\eean

ii) Nonergodic, zero entropy case. 

In this case $\ln\rho_{F}=0$. Thus if $F$ is nondiagonalizable then 
(\ref{na}) reads
\bean
\|P_{\vep,\alpha}^{n}\|\approx e^{-\vep R_{diss}^{-1} n} \|F^{n}\|\approx 
e^{-\vep R_{diss}^{-1}n + c_{1}\ln n +c_{2}},
\eean
with $0<c_{1}\leq d$.
This immediately yields
\bean
n_{p}\approx  R_{diss}\frac{c_{1}}{\vep},\qquad \frac{n_{th}}{\ln(n_{th})} \sim \frac{1}{\vep}.
\eean
And moreover $\|P_{\vep,\alpha}^{n_{p}}\|\sim (1/\vep)^{c_{1}}$.

If $F$ is diagonalizable then $\|F^{n}\|=1$ and in
this case $\|P_{\vep,\alpha}^{n}\|\approx e^{-\vep R_{diss}^{-1} n}$ which implies 
$n_{p}=0$. 

iii) If $F$ is diagonalizable, then form (\ref{ulestimateA}) we know that for 
any $0<\del<1$ and sufficiently large $n$
\bea
\label{est}
\lam_{-\del}^{n}=e^{(1-\del)2\alpha\hat{h}(A) n} \leq
\min_{0 \not=\bk \in \IZ^{d}} \sum_{l=1}^{n}|A^{l}\bk|^{2\alpha} \leq
e^{(1+\del)2\alpha \hat{h}(A)n}=\lam_{+\del}^{n}.
\eea
Thus for large $n$ we have
\bean
\max_{n} e^{-\vep \lam_{+\del}^{n}}\rho_{F}^{n} \leq \max_{n}\|P_{\vep,\alpha}^{n}\|
\leq \max_{n} e^{-\vep \lam_{-\del}^{n}}\rho_{F}^{n}.
\eean
We obtain the following constraints for  $n_{p}$
\bean
\frac{1}{\ln \lam_{+\del}}\ln\left(\frac{\ln \rho_{F}}{\ln \lam_{+\del}}\right)
+\frac{1}{\ln \lam_{-\del}}\ln\left(\frac{1}{\vep}\right)
\leq n_{p}\leq \frac{1}{ln \lam_{-\del}}\ln\left(\frac{\ln \rho_{F}}{\ln \lam_{-\del}}\right)
+\frac{1}{\ln \lam_{-\del}}\ln\left(\frac{1}{\vep}\right).
\eean
This gives
\bean
\frac{1}{\ln \lam_{+\del}} \leq \lim_{\vep \rightarrow 0} \frac{n_{p}}{\ln(1/\vep)}
\leq\frac{1}{\ln \lam_{-\del}}.
\eean
Now since $\lam_{\pm \del} \rightarrow e^{2\alpha\hat{h}(F)}$ for $\del\rightarrow 0$ 
the above estimation yields the following asymptotics
\bean
n_{p} \approx \frac{1}{2\alpha\hat{h}(F)}\ln(1/\vep)\approx n_{diss}
,\qquad n_{th}-R_{diss}\ln(n_{th})\approx n_{diss}.
\eean
Similar asymptotic estimates (except for the constant) hold
for nondiagonalizable $F$.
$\qquad \blacksquare$

\bigskip
\appendix
\centerline{\Large \bf Appendices}
\section{Affine transformations}
\label{aAT}
In this appendix we present a slight generalization of the results
obtained in the paper. We consider here general affine transformations
of the torus.
The term {\em affine transformations} will be used here to refer to  
homeomorphisms of the torus with zero periodic 
but not necessary zero constant part (cf. Section \ref{sDTA}) i.e. 
transformations of the form $\tilde{F}=F+\bc$, where $F$ is a toral 
automorphism and $\bc$ is a constant shift vector.
 
We begin with a short discussion of the ergodicity of affine
transforms.
 
The relation between ergodicity of a given affine
transform $\tilde{F}$ and associated with it toral 
automorphism $F$ is summarized in
the following proposition (for the proof we refer to appendix B)
\begin{prop}
\label{affinerg}
Let $F$ be any toral automorphism. Then

i) If $F$ is ergodic then $\tilde{F}$ is also ergodic.

ii) If $F$ is not ergodic then $\tilde{F}$ is ergodic iff $1$ is the only
root of unity in the spectrum of $F$ and $\bc\cdot\bk\not\in\IZ^{d}$ 
for any integer eigenvector $\bk$ of $F^{\dagger}$.
\end{prop}
{\bf Proof.} i) Assume $F$ is ergodic and for some $\bc$, $\tilde{F}=F+\bc$ is not ergodic.
Then there exists non-constant $f\in L_{0}^{2}(\IT^{d})$ satisfying 
$f=f\circ \tilde{F}$ or in the Fourier representation 
\bea
\label{inv}
\sum_{\bk \in \IZ^{d}}\hat{f}(\bk)\bfe_{\bk}=\sum_{\bk \in \IZ^{d}}
 e^{2 \pi i A^{-1}\bk\cdot \bc}\hat{f}(A^{-1}\bk)\bfe_{\bk}
\eea
where $A=F^{\dagger}$. 
Comparing the absolute values of the coefficients we get 
\bea
\label{inv1}
|\hat{f}(\bk)|=|\hat{f}(A^{-n}\bk)|
\eea
for any integer $n$ and any $\bk$.
However, ergodicity of $F$ implies that $A^{-n}\bk\not=\bk$ for all $\bk\not=0$, which contradicts our assumption that $f\in L_{0}^{2}(\IT^{d})$.

ii) We will use the following fact, which can be proved by
simple application of rational canonical decomposition. 
For any $A\in SL(d,\IZ)$ the following conditions are equivalent

a) $A$ possesses in its spectrum a root of unity not equal to one.

b) There exists nonzero $\bk \in \IZ^{d}$ and a positive integer $n$ such
that $\bk+A\bk+...+A^{n-1}\bk=0$.
 
Now assume that $1$ is the only root of unity in spectrum of $F$ (and hence
of $A$) and $\bc\cdot\bk\not\in\IZ^{d}$ for any integer eigenvector $\bk$ of 
$A$, and that both $F$ and $\tilde{F}$ are not ergodic. 
The latter assumption implies the 
existence of a non-constant $f\in L_{0}^{2}(\IT^{d})$ satisfying 
equations (\ref{inv}) and (\ref{inv1}). 
Relation (\ref{inv1}) clearly implies that if
$\hat{f}(\bk)\not= 0$ then $A^{n}\bk=\bk$ for some $n$. Moreover, 
since $1$ is the only root of unity in spectrum of $A$,
we have, in view of b) that $A\bk=\bk$. Thus the only possible non-constant
invariant functions of $\tilde{F}$ are single Fourier modes 
$\bfe_{\bk}$ corresponding to integer eigenvectors of $A$. 
But if such a Fourier mode is invariant under $\tilde{F}$ then
directly form (\ref{inv}) one concludes that $e^{2 \pi i \bk\cdot \bc}=1$
or equivalently $\bk\cdot \bc \in\IZ^{d}$, for some integer eigenvector of $A$. 
To prove the converse we assume that $F$ is not ergodic and consider two
cases:

Case 1. $A$ possesses in its spectrum a root of unity not equal to one.
In this case according to condition b) there exists nonzero 
$\bk \in \IZ^{d}$ and a positive integer $n$ such
that $\bk+A\bk+...+A^{n-1}\bk=0$, which implies in particular that 
$A^{n}\bk=\bk$ and $A\bk\not=\bk$. Now we define the function
\bean
f=\bfe_{\bk}+e^{2\pi i \bk \cdot \bc}\bfe_{A\bk}+...+
e^{2\pi i \left(\sum_{l=0}^{n-2}A^{l}\bk \right)\cdot \bc}\bfe_{A^{n-1}\bk}
\eean
which clearly satisfies the condition $f=f\circ \tilde{F}$.
This proves that $\tilde{F}$ is not ergodic.

Case 2. There exists integer eigenvector of $A$ such that $\bk\cdot \bc \in\IZ^{d}$. Then clearly for such $\bk$, $f=e_{\bk}$ is $\tilde{F}$-invariant and
hence again $\tilde{F}$ is not ergodic. 
$\qquad \blacksquare$

We recall that $\mf{c}=(c_{1},..,c_{d})$ generates ergodic shift
on the torus iff $1,c_{1},..,c_{d}$ are linearly independent over rationals.
Thus as a direct consequence of the above proposition we get
\begin{cor}
\label{corerg}
If $F$ is not ergodic and $1$ is the only
root of unity in the spectrum of $F$ then $\tilde{F}$ is ergodic for all
vectors $\bc$ generating ergodic shifts on the torus.
\end{cor}

Now we are in a position to state and prove the generalization of 
Theorem \ref{thm2} from Section \ref{sDTA} to the case of affine transforms
(the corresponding generalizations of Theorem \ref{thm3} and 
Corollary \ref{cat} are straightforward)  

\begin{thm}
\label{thm2'}

Let $\tilde{F}$ be any affine transformation on the torus $\IT^{d}$, 
$F$ associated with $\tilde{F}$ toral automorphism and 
$T_{\vep,\alpha}=G_{\vep,\alpha}U_{\tilde{F}}$. Then
 
i) $T_{\vep,\alpha}$ has simple dissipation time iff $F$ is not ergodic.

ii) $T_{\vep,\alpha}$ has logarithmic dissipation time iff $F$ is ergodic. 

iii) If $T_{\vep,\alpha}$ has logarithmic dissipation time then the dissipation rate constant 
     satisfies the following constraint
\bean
\frac{1}{2 \alpha \hat{h}(F)} \leq R_{diss} \leq \frac{1}{2 \alpha \tilde{h}(F)},
\eean
where $\tilde{h}(F)\leq\hat{h}(F)$ is certain positive constant. 
\end{thm}

\begin{rem}
\label{remerg}
The dissipation time of an affine transformation $\tilde{F}$ is determined 
by ergodic properties of its linear part $F$ and hence not by
ergodic properties of $\tilde{F}$ itself. In particular all ergodic affine
transformations associated with nonergodic toral automorphisms
(cf. Proposition \ref{affinerg}) have simple dissipation time.     
\end{rem}

\textbf{Proof of Theorem \ref{thm2'}}
Specializing the general calculations of dissipation time presented 
in Section \ref{sDCal} to the case of
affine transformations $\tilde{F}=F+\bc$, with nonzero $\bc$, 
one easily finds the following counterparts of formulas (\ref{linear}) 
and (\ref{linear1})
\bean
u_{\bk,\bk'}&=&e^{2\pi i \bk \cdot \bc}\del_{A\bk,\bk'}, \\
\ml{U}_{n}(\bk_{0},\bk_{n})
&=&e^{2\pi i \left(\sum_{l=0}^{n-1}A^{l}\bk\right) \cdot \bc}
e^{-\vep\sum_{l=1}^{n}|A^{l}\bk|^{2\alpha}}\del_{A^{n}\bk_{0},\bk_{n}}.
\eean 
Now, in order to determine the dissipation time of 
$T_{\vep,\alpha}=G_{\vep,\alpha}U_{\tilde{F}}$ one has to determine the 
asymptotics of $\|T_{\vep,\alpha}^{n}\|$ as $n$ goes to infinity.  
According to the above formulas and formulas (\ref{uppernorm}) and (\ref{Tnorm3})
from Section \ref{sDCal} the value of $\|T_{\vep,\alpha}^{n}\|$ does not depend on $\bc$,
which reduces the proof to the case $\bc=0$ considered in the main body
of the paper.  $\qquad \blacksquare$

\section{Proofs of some elementary facts}
\label{aP}
\textbf{Proof of Proposition \ref{pdiss}}

The proof will be based on the Riesz convexity theorem (see \cite{Zyg}, pp. 93-100) which
states that for any operator $T$ defined
on $L^p(\IT^d), 1\leq p\leq \infty,$  $\ln{\|T\|_p}$
is a convex function of $p^{-1}$. 
On the space $L^p(\IT^d)$ we consider the operator
 $\widetilde{T}:=T_{\vep,\alpha}-\la\cdot\ra$ and we have 
the relation $\widetilde{T}^n f= T^n_{\vep,\alpha} (f-\la f\ra), \forall f\in L^p(\IT^d), n\geq 1$ 
because $T_{\vep,\alpha}$ is conservative.
Now since $\|f -\la f\ra\|_p\leq 2\|f\|_p$, it follows
that
\bea
\label{14}
\|\widetilde{T}^n\|_p &\leq &2 \|T^n_{\vep,\alpha}\|_{p,0}\leq 2 \\
\label{15}
\|T^n_{\vep,\alpha}\|_{p,0} &\leq &\|\widetilde{T}^n\|_p
\eea
for $1\leq p\leq \infty, n\geq 1$.
The Riesz convexity theorem implies that if $p<q<\infty$
\bea
\label{16}
\ln{\|\widetilde{T}^{n}\|_q}\leq \frac{p}{q}\ln{\|\widetilde{T}^{n}\|_p}+\left(1-\frac{p}{q}
\right)\ln{\|\widetilde{T}^{n}\|_\infty}
\eea
while if $1<q<p$
\bea
\label{17}
\ln{\|\widetilde{T}^{n}\|_q}\leq \left(\frac{1-1/q}{1-1/p}\right)
\ln{\|\widetilde{T}^{n}\|_p}+\left(1-\frac{1-1/q}{1-1/p}
\right)\ln{\|\widetilde{T}^{n}\|_1}.
\eea
>From (\ref{16})-(\ref{17}) we have the interpolation relations
\bea
\label{16'}
\|\widetilde{T}^n\|_q &\leq& \|\widetilde{T}^n\|_p^{p/q}\|\widetilde{T}^n\|_\infty^{1-p/q}, \quad p<q<\infty\\
\label{17'}
\|\widetilde{T}^n\|_q &\leq &
\|\widetilde{T}^n\|_p^{(1-q^{-1})/(1-p^{-1})}\|\widetilde{T}^n\|_1^{1-(1-q^{-1})/(1-p^{-1})},\quad
1<q<p
\eea
which, along with (\ref{14})-(\ref{15}), imply
\bean
\|T^n_{\vep,\alpha}\|_{q,0}&\leq& 2\|T^n_{\vep,\alpha}\|_{p,0}^{p/q},\quad p<q<\infty\\
\|T^n_{\vep,\alpha}\|_{q,0}&\leq& 2\|T^n_{\vep,\alpha}\|_{p,0}^{(1-q^{-1})/(1-p^{-1})},
\quad 1<q<p
\eean
This proves that the order of divergence
of $n_{diss}(p)$ are the same for $1<p<\infty$.
Estimates (\ref{16'})-(\ref{17'}) also show
that the order of divergence of
$n_{diss}(1)$ and $n_{diss}(\infty)$ is at least as high as
$n_{diss}(p), 1<p<\infty$. $\qquad \blacksquare$

\bigskip
\textbf{Proof of Lemma \ref{gub}}

Using the notation introduced in Section \ref{sDCal} one has
\bean
&&T_{\vep,\alpha}^{n}\bfe_{\bk_{0}}=
(G_{\vep,\alpha}U)^{n}\bfe_{\bk_{0}}
=(G_{\vep,\alpha}U)^{n-1}\sum_{0 \not = \bk_{1} \in \IZ^{d}}u_{\bk_{0},\bk_{1}}
e^{-\vep |\bk_{1}|^{2\alpha}}\bfe_{\bk_{1}} \\ 
&=& \sum_{0 \not = \bk_{1},...,\bk_{n} \in \IZ^{d}}u_{\bk_{0},\bk_{1}}
u_{\bk_{1},\bk_{2}}...u_{\bk_{n-1},\bk_{n}} e^{-\vep \sum_{l=1}^{n}|\bk_{l}|^{2\alpha}}
\bfe_{\bk_{n}}=\sum_{0 \not = \bk_{n} \in \IZ^{d}} 
\ml{U}_{n}(\bk_0,\bk_{n})\bfe_{\bk_{n}}.
\eean

We note that for any $n$ and $\bk_{n} \in \IZ^{d}$, the sequence 
$\ml{U}_n(\bk_{0},\bk_{n})$ (indexed by $\bk_{0}\in \IZ^{d}$) belongs to $l^{2}(\IZ^{d})$.
Indeed, using the Cauchy-Schwarz inequality and identity (\ref{ul2}) one gets for $n=2$,
\bean
&&\sum_{0 \not = \bk_{0} \in \IZ^{d}}|\ml{U}_2(\bk_{0},\bk_{2})|^{2}=
\sum_{0 \not = \bk_{0} \in \IZ^{d}}\left|\sum_{0 \not = \bk_{1} \in \IZ^{d}}
u_{\bk_{0},\bk_{1}}u_{\bk_{1},\bk_{2}} 
e^{-\vep (|\bk_{1}|^{2\alpha}+|\bk_{2}|^{2\alpha})}\right|^{2} \\
&\leq&\sum_{0 \not = \bk_{0} \in \IZ^{d}}\sum_{0 \not = \bk_{1} \in \IZ^{d}}
|u_{\bk_{0},\bk_{1}}|^{2}e^{-\vep |\bk_{1}|^{2\alpha}}
\sum_{0 \not = \bk_{1} \in \IZ^{d}}|u_{\bk_{1},\bk_{2}}|^{2}
e^{-\vep |\bk_{1}|^{2\alpha}}e^{-2\vep |\bk_{2}|^{2\alpha}} \\
&\leq&\sum_{0 \not = \bk_{1} \in \IZ^{d}}e^{-\vep |\bk_{1}|^{2\alpha}}
\sum_{0 \not = \bk_{1} \in \IZ^{d}}e^{-\vep |\bk_{1}|^{2\alpha}}
e^{-2\vep |\bk_{2}|^{2\alpha}} = Ke^{-2\vep |\bk_{2}|^{2\alpha}},
\eean
where K denotes a constant. Similar estimates hold for $n>2$.

Now applying the Cauchy-Schwarz inequality in (\ref{Tnorm1}) we get
\bea
\|T_{\vep,\alpha}^{n}f\|^{2} &\leq&
\sum_{0 \not = \bk_{n} \in \IZ^{d}}
\sum_{\bk_{0} \in \ml{S}_{n}(\bk_{n})}
|\hat{f}(\bk_{0})|^{2}
\sum_{\bk_{0} \in \ml{S}_{n}(\bk_{n})}
|\ml{U}_n(\bk_{0},\bk_{n})|^{2}.
\eea $\qquad \blacksquare$

\bigskip
\textbf{Proof of part i) of Theorem \ref{thm2}.}

In view of theorem \ref{nonergodic} it suffices to construct an  
eigenfunction of $U_{F}$ which belongs to 
$L^{2}_{0}(\IT^{d}) \cap H^{2\alpha}(\IT^{d})$.
Directly from Proposition \ref{ergodicTA} one concludes that $F$, and
hence also $A$, possesses a root of unity in its spectrum. This means
that $A^{m}\bk_{0}=\bk_{0}$, for some $m$ and certain nonzero vector
$\bk_{0}$, which can be chosen to be an integer.
Now we define
\bean
f=\bfe_{\bk_{0}}+\bfe_{A\bk_{0}}+...+\bfe_{A^{m-1}\bk_{0}}
\eean
Obviously $f\in L^{2}_{0}(\IT^{d}) \cap H^{2\alpha}(\IT^{d})$, 
for any $\alpha$. To complete the proof it suffices to notice that
\bean
U_{F}f=\bfe_{A\bk_{0}}+\bfe_{A^{2}\bk_{0}}+...+\bfe_{A^{m}\bk_{0}}
=\bfe_{\bk_{0}}+\bfe_{A\bk_{0}}+...+\bfe_{A^{m-1}\bk_{0}}=f.
\eean  $\qquad \blacksquare$

\bigskip
\textbf{Proof of Proposition \ref{irred}.}

For the purposes of the proof we use the following abbreviation
\begin{itemize}
\item $PRS(\IR^{d})$ - proper rational subspace of $\IR^{d}$.

\item $PIS(A,\IR^{d})$ - proper $A$-invariant subspace of $\IR^{d}$.

\item $PRIS(A,\IR^{d})$ - proper rational $A$-invariant subspace of $\IR^{d}$.
\end{itemize}
\begin{description}
\item[a) $\Rightarrow $ b)]. Suppose there exists $PRIS(A,\IR^{d})$ $S_{1}$. Let $A_{1}$ be a matrix 
representing invariant rational block associated with $S_{1}$. Then $A_{1}$ is rational 
matrix and its characteristic polynomial $P_{1}$ belongs to $\IQ[x]$. 
Let $P$ denote the characteristic polynomial of $A$. Then $P=P_{1}P_{2}$ and since both 
$P,P_{1}\in \IQ[x]$ then also $P_{2}\in \IQ[x]$, which means $P$ and hence $A$ is not irreducible. 

\item[b) $\Rightarrow $ c)]. Assume there exists $PRS(\IR^{d})$ $S$ contained in $PIS(A,\IR^{d})$ $V$. 
Take any rational vector $\bq\in S$ and let $d_{0}=dimV$ then the set 
$\{\bq,A\bq,..,A^{d_{0}-1}\bq\}$ spans $PRIS(A,\IR^{d})$.

\item[c) $\Rightarrow$ d)]. Assume that for given $\bq$ and an arithmetic sequence $n_{1},...,n_{d}$, the set 
$S=\{A^{n_{1}}\bq,A^{n_{2}}\bq,...,A^{n_{d}}\bq\}$ does not form a basis. Since for some fixed
integer $r$, $n_{l}=n_{1}+(l-1)r$, we have $A^{n_{l}}\bq=(A^{r})^{l-1}A^{n_{1}}\bq=(A^{r})^{l-1}\hat{\bq}$, 
where $\hat{\bq}=A^{n_{1}}\bq$. Now consider the biggest subset 
$S_{0}=\{\hat{\bq},A^{r}\hat{\bq},(A^{r})^{2}\hat{\bq},...,(A^{r})^{d_{0}-1}\hat{\bq}\}$ such 
that $d_{0}<d$ and $S_{0}$ is linearly independent. Obviously $S_{0}$ spans a $PRIS(A^{r})$ which
is also a $PRIS(A)$.

\item[d) $\Rightarrow $ a)]. 
Suppose that characteristic polynomial $P$ of $A$ is not irreducible in $\IQ[x]$. 
Then $P=P_{1}P_{2}$, with $P_{1},P_{2}\in \IQ[x]$. From the Cayley-Hamilton theorem we get that
$0=P(A)=P_{1}(A)P_{2}(A)$. Hence for any nonzero rational vector $\bq$, either 1) $P_{2}(A)\bq=0$ 
or 2) $\hat{\bq}:=P_{2}(A)\bq \not = 0$ and $P_{1}(A)\hat{\bq}=0$. Since $\max\{deg(P_{1},P_{2})\}<d$,
there exists a nonzero rational vector $\tilde{\bq}$ (namely $\bq$ or $\hat{\bq}$) such that 
the set of iterations $\{\tilde{\bq},A\tilde{\bq},A^{2}\tilde{\bq},...,A^{d-1}\tilde{\bq}\}$ does 
not form a basis of $\IR^{d}$.

\item[e) $\Rightarrow $ f)]. Assume there exist nonzero $\bq \in \IQ^{d}$ orthogonal to certain 
$PIS(A,\IR^{d})$ $V$. Then for any $n$ and any $f\in V$, $\la (A^{\dagger})^{n}\bq, f \ra= \la \bq, A^{n}f \ra  = 0$ 
and hence $S=\{\bq,A^{\dagger}\bq,(A^{\dagger})^{2}\bq,...,(A^{\dagger})^{d-1}\bq\}$, cannot form a basis, which in view 
of equivalence a) $\Leftrightarrow$ d) implies reducibility of $A^{\dagger}$. 

\item[f) $\Rightarrow$ g)]. Suppose there exists $PIS(A,\IR^{d})$ $V$ contained in certain 
$PRS(\IR^{d})$ $S$. Since $S$ is rational, $S^{\perp}$ is also rational. Consider any rational 
vector $\bq \in S^{\perp}$, then $\la \bq,f \ra = 0$ for any $f \in V$.

\item[g) $\Rightarrow$ b)]. If there exists $PRIS(\IR^{d})$, then this subspace is $A$-invariant and 
contained in $PRS(\IR^{d})$ i.e
in itself. 
\end{description}

Now since b) is equivalent to a) it is enough to establish the equivalence
between a) and e) to complete the proof. But the latter equivalence is obvious
in view of the fact that $A$ and $A^{\dagger}$ have the same characteristic
polynomial.
$\qquad \blacksquare$

\bigskip
\textbf{Proof of Proposition \ref{zeroentropy}.}

Suppose $A$ is a toral automorphism of zero entropy. The latter property is 
equivalent
to the fact that all the eigenvalues of $A$ are of modulus $1$. Let $P_{A}$ be a
 characteristic
 polynomial of $A$. Consider any irreducible over $\IZ$ factor $P$ of polynomial
 $P_{A}$ and construct
a toral automorphism $B$ such that its characteristic polynomial is equal to $P$
 .
Obviously all the eigenvalues of $B$ are also the eigenvalues of $A$, and each 
 eigenvalue of
$A$ can be found among eigenvalues of some matrix $B$ of this type.
 Irreducibility of $P$ implies irreducibility and hence diagonalizability of $B$.

 Thus for any nonzero vector $\bk \in \IZ^{d}$ and any positive integer $n$ the
 following estimate holds
 $|B^{n}\bf{k}| \leq |\bf{k}|$,
which implies the existence (for each $\bf{k}$) of some integer $r$ such that $
 B^{r}\bf{k}=\bf{k}$.

The latter shows that all the eigenvalues of $B$ (and hence also of $A$) are 
roots
 of unity.
$\qquad \blacksquare$

\bigskip
\textbf{Proof of Proposition \ref{red}.}

We first show that irreducible polynomials $P\in\IQ[x]$ do not have repeated roots.
Indeed suppose $\lam$ is a root of $P$ of multiplicity greater that 1, then $\lam$ is
also a root of a derivative polynomial $P'\in\IQ[x]$. Since the minimal polynomial of $\lam$
must divide both $P$ and $P'$ and $deg(P')<deg(P)$ one immediately concludes that $P$ is not irreducible.
Now, suppose $A \in GL(d,\IQ)$ is completely decomposable over $\IQ$ and let (\ref{blockdiag}) be its 
block diagonal decomposition into irreducible blocks.
Each $P_{A_{j}}$, as a characteristic polynomial of $A_{j}$, is irreducible over $\IQ$ and hence does
not possesses repeated roots, which implies diagonalizability of each $A_{j}$ and hence of $A$.
On the other hand if $A$ is diagonalizable then its minimal polynomial does not possesses repeated 
roots, which implies that all characteristic polynomials associated with elementary divisors are 
(first powers of) irreducible polynomials. 
This implies irreducibility of each block in representation (\ref{blockdiag}). 
$\qquad \blacksquare$

\bigskip
\textbf{Proof of Proposition \ref{distinct}.}

Let $P_{A}$ be the characteristic polynomial of an irreducible matrix $A \in GL(d,\IQ)$.
Since $P_{A}$ is an irreducible element of $\IQ[x]$ it does not possesses repeated roots
(see the proof of Proposition \ref{red}).
$\qquad \blacksquare$

\commentout{
\bigskip
\textbf{Proof of Proposition \ref{underg}.}

We already know that ergodic automorphisms have positive entropy.
Assume $F$ is indecomposable and let $P^{m}$ be its characteristic
polynomial, where $P$ is irreducible over $Q$.
Let $B$ be any toral automorphism with characteristic polynomial $P$.
$B$ is irreducible and has positive entropy. Suppose $B$ is not 
ergodic. Take integral vector $k$ which satisfies
$B^{r}\bk=\bk$ for certain $r$. Then for any $n$ $|B^{n}\bk|<C$ which means
$\bk$ has no components in expanding eigenspace of $B$.
Therefore the set of consecutive iterates ${\bk,A\bk,A^{2}\bk,...}$ spans 
proper rational $B$-invariant subspace, which contradicts irreducibility
of $B$. Now ergodicity of $B$ clearly implies ergodicity of $F$. 
$\qquad \blacksquare$
}

\end{document}